\documentclass[11pt]{article}
\usepackage[a4paper,margin=27mm]{geometry}
\usepackage[T1]{fontenc}
\usepackage{lmodern}
\usepackage{amsmath,amssymb,amsthm,mathtools,booktabs,array}
\usepackage{microtype}
\usepackage{enumitem}
\usepackage{needspace}
\usepackage[colorlinks=true,linkcolor=blue,citecolor=blue,urlcolor=blue]{hyperref}
\newtheorem{theorem}{Theorem}[section]
\newtheorem{lemma}[theorem]{Lemma}
\newtheorem{proposition}[theorem]{Proposition}
\newtheorem{corollary}[theorem]{Corollary}
\newtheorem{claim}{Claim}
\newenvironment{claimproof}{%
  \begin{proof}[Proof of the claim]}%
  {\end{proof}}
\theoremstyle{definition}

\theoremstyle{remark}
\newtheorem{remark}[theorem]{Remark}
\DeclareMathOperator{\sat}{sat}
\DeclareMathOperator{\ex}{ex}
\DeclareMathOperator{\ES}{ES}
\DeclareMathOperator{\deficit}{def}
\newcommand{\SF}{\mathcal{S}}
\newcommand{\KF}{\mathcal{K}}
\newcommand{\KFive}{K_5^{(3)}}
\newcommand{\Low}{\mathcal{L}}
\newcommand{\High}{\mathcal{U}}
\newcommand{\lk}{L}
\newcommand{\sh}{\partial_2}
\newcommand{\sq}{\sqcup}
\setlist[enumerate]{label=\textup{(\arabic*)},leftmargin=*,itemsep=3pt}
\title{Saturation in three-uniform Berge hypergraphs:\\
the five-star spectrum and the \(K_4\) saturation number}
\author{XuanYu Cao\thanks{School of Mathematics and Statistics,
Central China Normal University, Wuhan, China.
E-mail: \texttt{2025112111@mails.ccnu.edu.cn}.}
\and Yuan Lai\thanks{School of Mathematics and Statistics, Central China Normal University, Wuhan 430079, Hubei, China. E-mail: \texttt{laiyuan0226@mails.ccnu.edu.cn}.}
\and Chunxiang Wang\thanks{School of Mathematics and Statistics,
and Hubei Key Lab--Math. Sci., Central China Normal University,
Wuhan, China. E-mail: \texttt{wcxiang@mail.ccnu.edu.cn}.}}
\date{}
\hypersetup{pdftitle={Saturation in three-uniform Berge hypergraphs: the five-star spectrum and the K4 saturation number},pdfauthor={XuanYu Cao, Yuan Lai, Chunxiang Wang}}
\begin{document}
\maketitle
\begin{abstract}
For a graph \(G\), a hypergraph is Berge-\(G\)-saturated if it contains no Berge copy of
\(G\), but adding any missing hyperedge creates one. We determine
the complete saturation spectrum of three-uniform Berge five-stars
for every order. For \(n\ge7\), every integer from
\(\lceil4n/3\rceil-3\) to \(2n-5\) occurs; the additional values
depend only on \(n\) modulo five. The structural ingredient is a
sharp bound on the remainder after deleting all complete five-vertex
components: for an \(n\)-vertex, \(m\)-edge Berge-five-star-saturated hypergraph,
this remainder has at most \(15(2n-m)/7\) vertices. We characterize
all equality cases. We also prove that every three-uniform
Berge-\(K_4\)-saturated hypergraph on \(n\ge7\) vertices has at
least \(n\) edges. Together with known constructions, this determines
its saturation number for every order and settles the conjecture of
Chen, He, and Xie; the new lower bound covers the previously unresolved
orders \(9\le n\le95\). Both arguments use the distinct representatives
forced by saturation to obtain local incidence constraints. For
five-stars these constraints bound the size of the remainder
and characterize equality. For \(K_4\), they reduce a putative
counterexample to a configuration on at most eight vertices; incidence
estimates and a sharp bound on the vertices covered by short cycles
in a bipartite graph exclude all remaining cases.

\end{abstract}

\section{Introduction}\label{sec:introduction}
We write \(K_t\) for the complete graph on \(t\) vertices and
\(K_{1,t}\) for the star with \(t\) leaves.
Saturation imposes a local condition on every missing edge, but its
extremal parameters describe the whole hypergraph. Exact results
therefore require more than a degree estimate: one must determine
how the configurations that complete different missing edges can
overlap. In this paper we use these overlaps to solve two problems
for three-uniform Berge hypergraphs: the complete saturation spectrum
of the five-star at every order, and the exact saturation number of
\(K_4\) at every order. The common step is to turn the representatives
required by saturation into incidence inequalities. The resulting
structural restrictions also explain the exceptional small cases.

Graph saturation asks how few edges a maximal forbidden-subgraph-free
graph can have. The classical work of Erd\H{o}s, Hajnal, and Moon~\cite{EHM}
determined the saturation number of complete graphs. K\'aszonyi and
Tuza~\cite{KT} subsequently proved a linear upper bound for the saturation
number of every fixed graph with at least one edge. For a broader account
of saturation problems and their variants, see the survey of Currie,
Faudree, Faudree, and Schmitt~\cite{CFFS}.

Saturation spectra for ordinary graphs have also been studied extensively.
Barefoot, Casey, Fisher,
Fraughnaugh, and Harary~\cite{BCFFH} studied the possible sizes of maximal
triangle-free graphs. Faudree, Faudree, Gould, Jacobson, and
Thomas~\cite{FFGJT} determined the saturation spectra of ordinary stars,
showing that every integer between the two endpoints occurs for
\(K_{1,t}\) when \(t\ge3\) and \(n\ge t+1\). The five-star Berge spectrum
considered here exhibits gaps near its upper endpoint.

The systematic study of Berge saturation was initiated by English,
Gordon, Graber, Methuku, and Sullivan~\cite{EGGMS}, who treated several
basic graph families, including stars and triangles. English, Gerbner,
Methuku, and Tait~\cite{EGMT} proved that the saturation number of Berge
copies of any fixed graph grows at most linearly when the uniformity is
between three and five. Gerbner, Patk\'os, Tuza, and Vizer~\cite{GPTV}
established linear bounds for further graph classes at arbitrary
uniformity. Xie~\cite{Xie} proves a linear upper
bound for every fixed uniformity and every fixed finite simple forbidden
graph with at least one edge. These results provide the general growth
scale; exact values and the possible sizes of saturated hypergraphs
require more detailed structural information.

\subsection{Berge saturation and notation}
A finite simple three-uniform hypergraph, or \emph{3-graph}, is a
hypergraph \(H\) with \(E(H)\subseteq\binom{V(H)}3\), where
\(\binom{X}{r}\) denotes the family of \(r\)-element subsets of \(X\).
We write \(n=|V(H)|\) and \(m=|E(H)|\) for its order and size.
A \emph{missing triple} is a member of
\(\binom{V(H)}3\setminus E(H)\).
A 3-graph \(H\) contains a \emph{Berge copy} of a graph \(G\) if there are an injection
\(\iota:V(G)\to V(H)\) and an injection \(\phi:E(G)\to E(H)\) such that
\(\{\iota(u),\iota(v)\}\subseteq\phi(uv)\) for every \(uv\in E(G)\).
We call \(\iota(V(G))\) the \emph{core vertex set} of the copy.
The hyperedge \(\phi(uv)\) is the \emph{representative} of the
\emph{core pair} \(\{\iota(u),\iota(v)\}\).
A \emph{Berge star} is a Berge copy of \(K_{1,t}\); the images
of its centre and leaves are called its centre and leaves.
A \emph{Berge triangle} is a Berge copy of \(K_3\), with three
distinct representative hyperedges.
The family of three-uniform Berge copies of \(G\) is denoted
by \(\operatorname{Berge}(G)\).
This extension of Berge paths and cycles to arbitrary core graphs was
developed by Gerbner and Palmer~\cite{GP}.
A 3-graph is \emph{\(\operatorname{Berge}(G)\)-saturated} if it contains no
Berge copy of \(G\), but adding any missing triple on its vertex set creates one.

For a family \(\mathcal F\) of 3-graphs, an \(\mathcal F\)-free
3-graph contains no member of \(\mathcal F\) as a subhypergraph.
It is \(\mathcal F\)-saturated if adding any missing triple
destroys this property. Its three-uniform saturation spectrum is
\[
 \ES_3(n,\mathcal F)=
 \{|E(H)|:|V(H)|=n,\ H\text{ is an }\mathcal F\text{-saturated 3-graph}\}.
\]
Its minimum is the saturation number \(\sat_3(n,\mathcal F)\).
Its maximum is the Tur\'an number \(\ex_3(n,\mathcal F)\), since every
maximum-size \(\mathcal F\)-free 3-graph is saturated. Determining the
spectrum therefore refines both endpoint problems.

Throughout the paper the two forbidden families are fixed as \(\SF=\operatorname{Berge}(K_{1,5}), \; \KF=\operatorname{Berge}(K_4)\). Here \(K_{1,5}\) is the star with five leaves and \(K_4\) is the
complete graph on four vertices.

The degree \(d_H(v)\) counts edges containing \(v\), and the
codegree \(d_H(u,v)\) counts edges containing \(\{u,v\}\).
We write \(\delta(H)\) for the minimum degree and call a
degree-three vertex \emph{cubic}. Subscripts are omitted when
the ambient hypergraph is fixed.

The \emph{two-shadow} \(\sh H\) is the ordinary graph on \(V(H)\)
whose edges are pairs contained in hyperedges of \(H\).
Neighbourhoods, adjacency, and components of \(H\) refer to this
graph: \(N_H(v)\) is the open neighbourhood and
\(N_H[v]=N_H(v)\cup\{v\}\).
The \emph{link graph} \(\lk_H(v)\) has vertex set \(N_H(v)\)
and edge set \(\{xy:vxy\in E(H)\}\).
Thus a link has no isolated vertices; all links are taken in the
whole ambient hypergraph unless an induced subhypergraph is specified.
We write \(H[X]\) for the subhypergraph induced by \(X\), and
\(H+e\) for the addition of a missing triple.
For a graph or hypergraph \(G\), the notation \(G-e\) means
deletion of the edge \(e\), retaining the vertex set.

The notation \(uv\) means \(\{u,v\}\), and \(uvw\) means \(\{u,v,w\}\); strings of digits
are used only for single-digit labels.
The complete 3-graph on \(t\) vertices is \(K_t^{(3)}\), also
abbreviated to \(C_t\) in the explicit constructions.
We use \(\sq\) for vertex-disjoint union and \(qG\) for the
vertex-disjoint union of \(q\) copies of \(G\).
We allow the empty hypergraph with no vertices, set \(\binom n3=0\) for
\(0\le n<3\), and interpret integer intervals as their
intersections with \(\mathbb Z\).

The \emph{incidence graph} of a hypergraph is the bipartite graph
whose two classes are its vertices and its edge identities,
with adjacency given by containment.
The \emph{cycle space} of an ordinary graph is the vector space
over the two-element field \(\mathbb F_2\) formed by edge sets
in which every vertex has even degree, with symmetric difference
as addition. Its dimension is the \emph{cycle rank}.
A graph is \emph{cyclic} if it contains a cycle and \emph{unicyclic}
if it is connected and contains exactly one cycle.

\subsection{Results and their relation to earlier work}
Austhof and English~\cite[Theorem~4]{AE} determined the saturation
numbers of Berge stars for sufficiently large orders. Their formula
specializes to \(\sat_3(n,\SF)=\lceil4n/3\rceil-3\), and their
lower-bound argument applies without a large-order assumption.
Gerbner, Methuku, and Palmer~\cite[Theorem~12]{GMP} proved
\(\ex_3(n,\SF)\le2n\), with equality when \(5\mid n\).
Bushaw, English, Heath, Johnston, and Rombach
\cite[Theorems~3.1 and~5.13, Proposition~5.1]{BEEJR} determined
the five-star saturation spectrum for sufficiently large orders divisible
by five, obtaining the interval \([\lceil4n/3\rceil-3,2n-5]\)
together with the isolated value \(2n\).
Our first result gives the spectrum for every order and every residue class.

\begin{theorem}\label{s-thm:spectrum}
For \(n\geq7\), put \(s_n=\lceil4n/3\rceil-3\). Then
\[
\ES_3(n,\SF)=
\begin{cases}
 \bigl([s_n,2n-5]\cap\mathbb Z\bigr)\cup\{2n\},&n\equiv0\pmod5,\\
 \bigl([s_n,2n-4]\cap\mathbb Z\bigr)\cup\{2n-2\},&n\equiv1\pmod5,\\
 [s_n,2n-4]\cap\mathbb Z,&n\equiv2,4\pmod5,\\
 [s_n,2n-5]\cap\mathbb Z,&n\equiv3\pmod5.
\end{cases}
\]
For \(n=6\), the spectrum is \(\{6,7,8,10\}\). For \(0\leq n\leq5\), it consists of the single value \(\binom n3\).
\end{theorem}

\begin{corollary}\label{s-thm:saturation}
For every nonnegative integer \(n\),
\[
 \sat_3(n,\SF)=
 \begin{cases}
  \binom n3,&0\leq n\leq5,\\
  6,&n=6,\\
  \left\lceil4n/3\right\rceil-3,&n\geq7.
 \end{cases}
\]
\end{corollary}

\begin{corollary}\label{s-thm:extremal}
For every nonnegative integer \(n=5q+t\), where \(0\leq t\leq4\),
we have \(\ex_3(n,\SF)=10q+\binom t3\).
The disjoint union \(qK_5^{(3)}\sqcup K_t^{(3)}\) attains this value.
\end{corollary}

Theorem~\ref{s-thm:spectrum} removes both the divisibility condition
and the large-order assumption in the previously known five-star
spectrum. Its upper endpoints give the exact Tur\'an remainders in
Corollary~\ref{s-thm:extremal}. The large-order saturation formula
itself is due to Austhof and English~\cite{AE}.

The structural contribution is separate from the interval constructions.
For an \(\SF\)-saturated 3-graph \(H\), define the \emph{deficit}
by \(\deficit(H)=2n-m\), and let the
\emph{reduced core} \(F(H)\) be obtained by deleting all components
isomorphic to \(K_5^{(3)}\), retaining every other component and every
isolated vertex. This reduced core is a hypergraph, distinct from
the core vertex set of a Berge copy. If \(r\) is the number of degree-six vertices in
\(F(H)\), Theorem~\ref{s-thm:core} proves the stronger inequality \(7|V(F(H))|+5r\le15\deficit(H)\). In particular, \(|V(F(H))|\le15\deficit(H)/7\). We determine
all equality cases in this latter bound: apart from the empty core,
they are built from disjoint seven-edge configurations, each with
two designated vertices. Additional edges partition the designated
vertices into triples. These systems are described precisely in
Subsection~\ref{s-sec:equality}. A construction from~\cite{BEEJR}
already attains the coefficient \(15/7\); the new statements are the
upper bound and its full equality classification. The link method
and the constructions from~\cite{BEEJR} are identified where used.

For Berge cliques, English, Kritschgau, Nahvi, and Sprangel~\cite{EKNS}
determined asymptotic saturation numbers and asked for the exact value
for three-uniform Berge-\(K_4\). Chen, He, and Xie~\cite{CHX} proved
\(\sat_3(n,\KF)=n\) for \(n=5,7,8\) and \(n\ge96\), determined
the exception at \(n=6\), and conjectured equality for every \(n\ge7\).
Their constructions already give the upper bound at all these orders.
We prove the matching lower bound for every \(n\ge7\), thereby
settling the remaining orders \(9\le n\le95\).

\begin{theorem}\label{k-thm:main}
For every nonnegative integer \(n\),
\[
 \sat_3(n,\KF)=
 \begin{cases}
  \binom n3,&0\le n\le4,\\
  5,&n=5,6,\\
  n,&n\ge7.
 \end{cases}
\]
\end{theorem}

The proof retains the minimum-degree reduction and local
framework of~\cite{CHX}. The additional ingredients are incidence
estimates for degree-two vertices and a bound on short-cycle coverage
in bipartite graphs. These give a uniform analytic proof of the lower
bound, including the previously known orders seven and eight.
Thus the extension to all orders follows from structural restrictions,
rather than from an enumeration of the remaining hypergraphs.

\subsection{A common proof strategy}
In either problem, adding a missing triple creates a forbidden Berge
copy. Deleting the new representative leaves four old representatives
for a five-star, or five for \(K_4\) with one core pair removed.
The corresponding partial configurations are defined precisely in
Subsection~\ref{s-sec:preliminaries} and Section~\ref{k-sec:local}.
The requirement that all representatives be distinct is the source
of the local restrictions in both arguments.

For five-stars, assigning distinct endpoints to distinct edges of a
link determines the largest possible star at its centre.
Vertices of degrees five and six then have tightly constrained links.
After the complete five-vertex components are deleted, the remaining
dense configurations require enough neighbours of degree at most four
to force the reduced-core bound. Define the \emph{degree deficiency}
of a vertex \(v\) to be \(6-d_H(v)\).
The global count is \(\sum_{v\in V(H)}(6-d_H(v))=3\deficit(H)\).
Small deficits therefore leave only small reduced cores, whose
possible structures determine the gaps at the top of the spectrum.
Explicit starting hypergraphs and two saturation-preserving
disjoint-union operations then realize every remaining size.
These operations increase order and size by \((6,8)\) and \((5,10)\).

For \(K_4\), each missing triple contains a pair that can be completed
to a Berge \(K_4\) using five old hyperedges and the new triple.
If an endpoint of that pair has degree two, both of its incident
hyperedges must represent pairs joining it to the two other core
vertices. This forces incidences with neighbours of degree at least
three. We first exclude edges containing two degree-two vertices in
a putative counterexample. We then find a pair of degree-two vertices
for which no such five-edge partial copy exists. Their combined set
of neighbours of degree at least three has at most eight vertices.
Missing triples through this pair force Berge triangles on that set.

Replacing the remaining hyperedges by their intersections with this
set, while retaining their identities, turns these triangles into
six-cycles in an incidence graph. Incidence inequalities and bounds
on the vertices covered by short cycles exclude the remaining
configurations. Theorem~\ref{k-thm:short-cycle-support} states the
cycle bound, including its sharpness; the final proof lists all
parameter reductions and cases. Both parts thus proceed from local
completion constraints to a global incidence bound and then to a
finite structural analysis.

\subsection{Organization}
Section~\ref{s-sec:core-proof} develops the five-star structural bound,
and Section~\ref{s-sec:spectrum-proof} proves the full spectrum.
Section~\ref{k-sec:local} introduces the local tools for Berge-\(K_4\).
Section~\ref{k-sec:lower} develops the incidence analysis,
and Section~\ref{k-sec:constructions} completes the proof of
Theorem~\ref{k-thm:main}, including the known upper bounds and all
remaining lower-bound cases.

\section{The structure of Berge-five-star-saturated hypergraphs}\label{s-sec:core-proof}
In this section and the next, saturation and freeness refer to \(\SF\).
We use the deficit \(\deficit(H)=2|V(H)|-|E(H)|\) and the reduced core
\(F(H)\) defined in the introduction, and write \(d=\deficit(H)\) when
the ambient hypergraph is fixed.

\subsection{Links and elementary saturation properties}\label{s-sec:preliminaries}
For a finite simple graph \(L\), its \emph{bipartite incidence
graph} has vertex classes \(E(L)\) and \(V(L)\), with an edge
between \(e\in E(L)\) and \(u\in V(L)\) precisely when
\(u\in e\). The \emph{incidence-matching rank} of \(L\),
denoted by \(\rho(L)\), is the maximum size of a matching
in this bipartite graph. Equivalently, \(\rho(L)\) is the
largest integer \(r\) for which there are distinct edges
\(e_1,\ldots,e_r\in E(L)\) and distinct vertices
\(u_1,\ldots,u_r\in V(L)\) satisfying \(u_i\in e_i\)
for every \(i\).

For a vertex \(v\) of a 3-graph \(H\), the \emph{link rank}
at \(v\) means \(\rho(\lk_H(v))\). The \emph{Berge degree}
\(b_H(v)\) is the largest number of leaves in a Berge star
centred at \(v\). An incidence matching in \(\lk_H(v)\)
selects distinct hyperedges containing \(v\) and assigns
a distinct leaf to each of them. Consequently,
\(b_H(v)=\rho(\lk_H(v))\).

Let \(H\) be \(\SF\)-free, and let \(e\) be a missing triple
containing \(v\). A \emph{witness of Berge five-star} at \(v\)
for the addition of \(e\) consists of five distinct vertices
\(u_1,\ldots,u_5\in V(H)\setminus\{v\}\) and four distinct
hyperedges \(f_1,\ldots,f_4\in E(H)\) such that
\(
 \{v,u_i\}\subseteq f_i\;(1\le i\le4),
 \; \{v,u_5\}\subseteq e.
\)
Thus \(f_1,\ldots,f_4,e\) represent the five pairs of a Berge
five-star. The leaves need not have degree one in \(H\), and
the representative hyperedges may have additional intersections.
Since \(H\) is \(\SF\)-free, every Berge five-star in \(H+e\)
uses \(e\), and its centre belongs to \(e\).

A vertex \(v\in V(H)\) is \emph{universally sensitive} if,
for every missing triple \(e\) containing \(v\), there is a
witness of Berge five-star at \(v\) for the addition of \(e\).
Equivalently, adding any such triple creates a Berge five-star
centred at \(v\). An \(\SF\)-saturated 3-graph is
\emph{attachment-ready} if all its vertices have Berge degree four.
We collect the elementary facts needed throughout the proof in one lemma.

\begin{lemma}\label{s-lem:rank}\label{s-lem:sensitive}\label{s-lem:local}\label{s-lem:add}\label{s-lem:sat-lower}
The following facts concern incidence matchings and \(\SF\)-saturation.
\begin{enumerate}
\item If \(L\) has no isolated vertices, then
\(\rho(L)=|V(L)|-t(L)\), where \(t(L)\) counts its tree
components. A 3-graph is \(\SF\)-free precisely when every vertex
has link rank at most four.

\item A four-edge forest, \(K_3\sq K_2\), or \(K_4\) as a link
makes its vertex universally sensitive. At a vertex of Berge
degree four, adding an incident triple that introduces a new
neighbour creates a Berge five-star centred at that vertex.
Consequently, the disjoint union of an attachment-ready
\(\SF\)-saturated 3-graph and any \(\SF\)-saturated 3-graph
is \(\SF\)-saturated.

\item In an \(\SF\)-free 3-graph, every vertex has degree at most
six; a degree-six vertex has link \(K_4\), and a degree-five
vertex has link \(K_4-e\). In an \(\SF\)-saturated 3-graph,
two adjacent degree-six vertices belong to a \(\KFive\)
component.

\item Deleting any collection of \(\KFive\) components from
an \(\SF\)-saturated 3-graph preserves saturation and deficit.
Conversely, if \(F\) is \(\SF\)-saturated, then
\(F\sq q\KFive\) is \(\SF\)-saturated for every integer \(q\ge0\).

\item Every \(n\)-vertex \(\SF\)-saturated 3-graph has at least
\(\lceil4n/3\rceil-3\) edges. At order six it has at least
six edges.
\end{enumerate}
\end{lemma}

\begin{proof}
For (1), we prove the rank formula componentwise; see
also~\cite[Lemma~2.1]{BEEJR}. A tree component on \(t\)
vertices has rank at most \(t-1\), since it has \(t-1\)
edges. Root the tree and assign each edge to its endpoint
farther from the root. This gives an incidence matching
of size \(t-1\).

A connected component containing a cycle has rank at most
its number \(t\) of vertices. Choose a spanning unicyclic
subgraph, orient its cycle cyclically and its attached
trees away from the cycle, and assign each edge to its
head. Every vertex receives exactly one edge, giving an
incidence matching of size \(t\). Ranks add over components,
which proves the formula.

An incidence matching in \(\lk_H(v)\) corresponds precisely
to distinct representative hyperedges with distinct leaves
at \(v\). Thus \(b_H(v)=\rho(\lk_H(v))\), proving the
freeness criterion.

For (2), first suppose that \(b_H(v)=4\) and that a missing
triple \(e=vxy\) introduces a new neighbour, say
\(y\notin N_H(v)\). Retain an incidence matching of size
four in \(\lk_H(v)\) and assign the new link edge \(xy\)
to \(y\). Together with the new hyperedge, these four
representatives form a Berge five-star centred at \(v\).

Each of the three stated links has rank four. If the link
is a four-edge forest, adding any missing link edge gives
a five-edge graph in which every component is either a
tree or a unicyclic graph. By (1), its rank equals its
number of edges and is therefore five. If the link is
\(K_3\sq K_2\), a missing edge whose endpoints are both
already in the link must join its two components. The
result is a connected unicyclic graph on five vertices,
again of rank five. Otherwise, the new edge introduces
a new neighbour, and the preceding argument applies.
Finally, if the link is \(K_4\), every missing link edge
must introduce a new neighbour. This proves universal
sensitivity in all three cases; see
also~\cite[Lemma~2.3]{BEEJR}.

Now let \(G\) be attachment-ready and let \(H\) be any
\(\SF\)-saturated 3-graph. Their disjoint union is
\(\SF\)-free, since all representative hyperedges of a
Berge star contain its centre and hence lie in one
component. An internal missing triple is handled by
saturation of its factor. A crossing triple introduces
a new neighbour at any of its vertices in \(G\), where
the Berge degree is four. Hence \(G\sq H\) is saturated.

For (3), any family of at most four distinct edges of a
simple graph has at least as many endpoints as edges.
Thus every such family satisfies Hall's condition and
admits distinct endpoint representatives. For five
distinct edges, the only possible failure of Hall's
condition is that all five edges lie on four vertices.
In that case they form a \(K_4-e\).

Consequently, a link of rank at most four with at least
five edges contains a \(K_4-e\) on some four-set \(S\).
This subgraph has an incidence matching covering \(S\).
Any further edge with an endpoint outside \(S\) would
extend that matching to size five. All link edges
therefore lie on \(S\), so the link is either \(K_4-e\)
or \(K_4\). This proves the degree bound and the two
link descriptions.

Suppose now that \(H\) is saturated and that \(u,v\)
are adjacent degree-six vertices. Put \(T=N_H[u]\).
Since \(\lk_H(u)=K_4\), we have \(uvw\in E(H)\) for
every \(w\in T\setminus\{u,v\}\). Since
\(\lk_H(v)=K_4\) as well, it follows that \(N_H[v]=T\).

The hyperedges incident with \(u\) or \(v\) supply nine
of the ten triples of \(T\). Hence every vertex of \(T\)
has degree at least five in \(H[T]\) and four neighbours in \(T\). An edge meeting both \(T\) and its complement
would give some vertex of \(T\) degree at least six and
at least five neighbours, contrary to the degree bound
and the description of a degree-six link.

If the remaining triple of \(T\) were missing, adding it
would leave \(T\) as a five-vertex component and hence
could not create a Berge five-star. Saturation therefore
forces that triple, so \(T\) is a \(\KFive\) component.
These restrictions also appear
in~\cite[Observations~5.2--5.3 and Claim~5.5]{BEEJR}.

For (4), let \(H'\) be obtained from an \(\SF\)-saturated
3-graph \(H\) by deleting the specified components.
Clearly, \(H'\) is \(\SF\)-free. Let \(e\) be a missing
triple on \(V(H')\). Saturation of \(H\) gives a Berge
five-star in \(H+e\), and this star must use \(e\).
Its centre therefore belongs to \(e\). Every old
representative hyperedge contains that centre and
consequently belongs to \(H'\), since only whole
components were deleted. The same star occurs in
\(H'+e\), proving that \(H'\) is saturated.

Each deleted \(\KFive\) component removes five vertices
and ten edges, so the deficit \(2n-m\) is unchanged.
Conversely, \(\KFive\) is saturated and every one of its
vertices has Berge degree four. Repeated application
of (2) therefore shows that \(F\sq q\KFive\) is saturated
whenever \(F\) is saturated; see
also~\cite[Observation~3.4]{BEEJR}.

For (5), we give the lower-bound argument
from~\cite[proof of Theorem~7]{AE}, which does not require
a large-order assumption. Let \(H\) have \(n\) vertices
and \(m\) edges, and put
\(A=\{v\in V(H):d_H(v)\le3\}\) and \(a=|A|\).
If a triple \(e\subseteq A\) were missing, every Berge
five-star created by adding \(e\) would have its centre
in \(e\). Each vertex of \(e\) would then have degree
at most four, which is impossible. Hence \(H[A]\)
is complete.

If \(a\ge5\), every vertex of \(A\) would have degree at
least \(\binom{a-1}{2}\ge6\), a contradiction. Thus
\(0\le a\le4\). Summing degrees gives
\(3m\ge4(n-a)+3\binom a3\).
For \(a=0,1,2,3,4\), the term
\(-4a+3\binom a3\) equals \(0,-4,-8,-9,-4\),
respectively. Therefore \(3m\ge4n-9\), and hence
\(m\ge\lceil4n/3\rceil-3\).

It remains to strengthen the bound when \(n=6\).
The general bound gives \(m\ge5\). Suppose that \(m=5\).
The inequality \(15\ge4(6-a)+3\binom a3\) excludes
\(a=0,1,2,4\), so \(a=3\), and equality holds throughout
the degree estimate. The three vertices of \(A\)
therefore each have degree one and form an isolated
triple, while the other three vertices each have degree
four. This is impossible, since a simple 3-graph on
three vertices has at most one edge. Thus \(m\ge6\).
\end{proof}

\subsection{The sun construction}\label{s-sec:sun}
Let \(S\) be the sun of~\cite[Construction~4.2]{BEEJR} at
\(\ell=5\), equivalently the complete tripartite 3-graph with three
parts of size two. It has \((n,m,d)=(6,8,4)\).
By~\cite[Lemma~4.3]{BEEJR}, it is saturated and every vertex has
Berge degree four. Thus \(S\) is attachment-ready, and
\(S\sq H\) is saturated whenever \(H\) is saturated.

\subsection{A sharp bound on the reduced core}
\begin{theorem}\label{s-thm:core}
Let \(H\) be \(\SF\)-saturated, let \(F=F(H)\), and let \(r\) be the number of degree-six vertices in \(F\). Then \(F\) is \(\SF\)-saturated, \(\deficit(F)=\deficit(H)\), and \(7|V(F)|+5r\leq15\deficit(H)\).
In particular, \(|V(F)|\leq\lfloor15\deficit(H)/7\rfloor\). The coefficient \(15/7\) is sharp. The equality \(7|V(F)|=15\deficit(H)\) holds precisely for the empty core and the systems described in Subsection~\ref{s-sec:equality}.
\end{theorem}

\begin{proof}
By Lemma~\ref{s-lem:add}(4), the reduced core \(F\) is saturated and
has the same deficit as \(H\). Partition its vertices as
\(
 R=\{v:d_F(v)=6\},\; A=\{v:d_F(v)=5\},\; Q=\{v:d_F(v)\leq4\}.
\)
Throughout this section and the next, vertices in \(A\) are
called \emph{high} and those in \(Q\) \emph{low}; these terms
refer to ordinary hypergraph degrees five and at most four.
The \emph{high shadow} is \((\sh F)[A]\).
The \emph{block} of one of its components consists of all hyperedges
meeting that component, together with their vertices.
Degrees within this block are called its \emph{internal degrees}.
By Lemma~\ref{s-lem:local}(3), every neighbour of a vertex in \(R\)
has codegree three with it, and every neighbour of a vertex in \(A\)
has codegree at least two with it. The following claims separate the local structural argument from
the final count. First we isolate the neighbourhoods of degree-six
vertices. We then show that the shadow induced by degree-five vertices
has only isolated vertices, edges, and triangles, and identify the two
nontrivial blocks. The concluding double count uses only the resulting
restrictions on high--low adjacencies.

\begin{claim}\label{s-lem:R}
The sets \(N(z)\), for \(z\in R\), are pairwise disjoint sets of four vertices in \(Q\). No vertex in their union is adjacent to a vertex of \(A\).
\end{claim}
\begin{claimproof}
By Lemma~\ref{s-lem:local}(3), \(R\) is independent in the shadow. Suppose \(z\in R\) is adjacent to \(x\in A\), and write \(N(z)=\{x,a,b,c\}\). There are all six \(z\)-edges determined by pairs of these neighbours. The link of \(x\) is \(K_4-e\), so its two further edges are two of \(xab,xac,xbc\), say the first two. In the five-set \(T=\{z,x,a,b,c\}\) the degrees in \(F[T]\) are at least \(6,5,5,4,4\). Every vertex already has the other four as neighbours. An edge meeting both \(T\) and its complement would give a vertex degree at least five and at least five neighbours, contradicting Lemma~\ref{s-lem:local}(3). Hence \(T\) is a component, and saturation makes it a deleted \(\KFive\), a contradiction. Thus \(N(z)\subseteq Q\).

If \(q\in N(z)\cap N(z')\) for distinct \(z,z'\in R\), the three \(zq\)-edges and the three \(z'q\)-edges are six distinct edges, since no edge contains both \(z\) and \(z'\). This gives \(d(q)\geq6\), a contradiction. Finally, if \(q\in N(z)\) were adjacent to \(a\in A\), its three \(zq\)-edges and at least two \(aq\)-edges would again be distinct, since \(za\) is not a shadow edge. This contradicts \(d(q)\leq4\).
\end{claimproof}

\begin{claim}\label{s-lem:two-low}
Every vertex in \(A\) has at least two neighbours in \(Q\).
\end{claim}
\begin{claimproof}
There are no \(A\)--\(R\) shadow edges by Claim~\ref{s-lem:R}. Suppose \(v\in A\) has at most one neighbour in \(Q\), and let \(T=N[v]\), a five-set. At least three of its four neighbours are in \(A\). In \(\lk(v)=K_4-e\), a neighbour of degree three in the link already has the other four members of \(T\) as neighbours and therefore has no neighbour outside \(T\) if it lies in \(A\).

Let \(u\in A\cap N(v)\) have degree two in \(\lk(v)\). Among its three presently known neighbours are \(v\) and at least one neighbour of \(v\) in \(A\) that has degree three in \(\lk(v)\). Both are confined to \(T\). If the fourth neighbour \(z\) of \(u\) were outside \(T\), then in \(\lk(u)\) the vertex \(z\) could be adjacent to at most one of the other three vertices. Thus \(\lk(u)\) would have at most \(3+1=4\) edges, a contradiction. All vertices of \(A\cap T\) are consequently confined to \(T\).

If every member of \(T\) lies in \(A\), \(T\) is a five-vertex component and saturation is impossible. Otherwise write \(T\cap Q=\{q\}\). The four other vertices have degree five in \(F[T]\). Relative to the complete 3-graph on \(T\), let \(f\) missing triples avoid \(q\) and \(g\) contain \(q\). Counting their incidences with the four degree-five vertices gives \(3f+2g=4\), whose only nonnegative integer solution is \(f=0,g=2\). Thus \(q\) has degree four in \(F[T]\); it has no further edge since \(q\in Q\). Again \(T\) is a noncomplete five-vertex component, contradicting saturation.
\end{claimproof}

\begin{claim}\label{s-lem:Acomponents}
Every component of \((\sh F)[A]\) is \(K_1\), \(K_2\), or \(K_3\).
\end{claim}
\begin{claimproof}
Claim~\ref{s-lem:two-low} gives maximum degree two. Suppose \(v\in A\) has nonadjacent neighbours \(a,b\in A\) and remaining neighbours \(x,y\in Q\). Its link is \(K_4-ab\), giving the five edges \(vax,vay,vbx,vby,vxy\).
The fourth neighbour \(z\) of \(a\), in addition to \(v,x,y\), is not adjacent to \(v\): it is distinct from \(b\), since \(ab\) is not a shadow edge. Therefore \(\lk(a)=K_4-vz\), and \(axy,axz,ayz\) are edges. The five distinct edges \(vax,vbx,vxy,axy,axz\) all contain \(x\), contradicting \(x\in Q\). Hence the two neighbours of every degree-two vertex of \((\sh F)[A]\) are adjacent, proving the assertion.
\end{claimproof}

\begin{claim}\label{s-lem:three-high}
Every vertex of \(Q\) is adjacent to at most three vertices of \(A\).
\end{claim}
\begin{claimproof}
Every such adjacency has codegree at least two. If \(q\) had at least four neighbours in \(A\), its degree at most four would force exactly four, each of codegree two, with all four \(q\)-edges containing two vertices of \(A\). Its link on those four vertices would be a two-regular simple graph, hence a four-cycle. This gives a four-cycle in \((\sh F)[A]\), contradicting Claim~\ref{s-lem:Acomponents}.
\end{claimproof}

\begin{claim}The two nontrivial high-vertex blocks are as follows.\label{s-lem:blocks}
\leavevmode
\begin{enumerate}
\item A triangle component \(\{a,b,c\}\) of \((\sh F)[A]\) has two common neighbours \(x,y\in Q\). Its incident edges are exactly \(abc,abx,aby,acx,acy,bcx,bcy\). Each of \(x,y\) has internal degree three and is adjacent to no other vertex of \(A\).
\item An edge component \(\{a,b\}\) has three common neighbours in \(Q\). Its incident edges consist of all three \(abq\) edges and two triples on \(a\) and pairs of the three low vertices, and similarly two on \(b\). The omitted low pairs are different. The three low vertices have internal degrees \(4,4,3\) and are adjacent to no other vertex of \(A\).
\end{enumerate}
\end{claim}
\begin{claimproof}
For (1), write \(N(a)=\{b,c,x,y\}\) and \(\lk(a)=K_4-e\). If \(e=xy\), the \(a\)-edges already give each of \(b,c\) all four other vertices of this five-set as neighbours. If \(e=bc\), the \(a\)-edges give both low neighbours to \(b,c\), while \(bc\) is a shadow edge by the assumed high triangle. In either case their low neighbours are \(x,y\). The remaining case is a missing high--low pair; by symmetry suppose \(e=cx\). The \(a\)-edges are \(abc,abx,aby,acy,axy\). Vertex \(b\) already has neighbours \(a,c,x,y\) and is confined to the five-set, whereas \(c\) already has neighbours \(a,b,y\). If the fourth neighbour \(z\) of \(c\) were outside the five-set, then \(az,bz\) would not be shadow edges, so the link on \(\{a,b,y,z\}\) could have at most \(ab,ay,by,yz\), only four edges. This contradicts \(d(c)=5\). Hence \(z=x\). Thus in every case \(a,b,c\) have the same two low neighbours.

Let \(t\in\{0,1\}\) count the triple \(abc\), and let
\(b_2,c_2\) count triples on this five-set containing respectively
two vertices of \(A\) and one of \(Q\), or one of \(A\) and two of \(Q\). Then \(3t+2b_2+c_2=15,\; b_2\leq6,\; c_2\leq3,\; b_2+2c_2\leq8\). The only possibility is \(t=1,b_2=6,c_2=0\). This seven-edge configuration is the \emph{triangle block}.

For (2), each endpoint of the edge component \(ab\) has three low neighbours. If \(d(a,b)=2\), write the two edges containing both vertices as \(abx,aby\).
Let \(z,w\) be their respective third low neighbours, which may coincide. The links force \(axy,axz,ayz\) and \(bxy,bxw,byw\), so \(x\) lies in the five distinct edges \(abx,axy,axz,bxy,bxw\), impossible. Hence \(d(a,b)=3\), and the three low neighbours coincide. Each high endpoint has two additional low-pair edges. If both omit the same pair, the third low vertex has degree five. Distinct omitted pairs give internal degrees \(4,4,3\).

In either block, a low vertex of internal degree at least three cannot be adjacent to an additional high vertex, since that adjacency supplies at least two further edges and high components are distinct.
\end{claimproof}

We now count the degree deficiencies. Put \(a=|A|\) and \(q=|Q|-4r\),
so \(q\) counts the low vertices outside the pairwise disjoint
neighbourhoods in Claim~\ref{s-lem:R}, and \(|V(F)|=5r+a+q\).
Claims~\ref{s-lem:R}, \ref{s-lem:two-low}, and
\ref{s-lem:three-high} show that every vertex in \(A\) has at least two
neighbours among these \(q\) low vertices, each of which has at most
three neighbours in \(A\). Double counting these adjacencies gives
\(2a\leq3q\). Since all \(4r\) neighbours of \(R\) have degree at most four, \(3\deficit(F)=\sum_{v\in V(F)}(6-d_F(v))\geq8r+a+2q\). Consequently \(15\deficit(F)-7|V(F)|\geq5r-2a+3q\geq5r\),
which proves the required bound. Its sharpness and the equality
classification are established in Subsection~\ref{s-sec:equality}.
\end{proof}

\subsection{Equality systems}\label{s-sec:equality}
Retaining the notation \(R,A,Q\) used above, we complete the equality statement of Theorem~\ref{s-thm:core} by identifying all equality cases.
Take \(t\) disjoint five-sets, each consisting of a high triple and two
vertices called \emph{ports}. On each five-set place the seven-edge
triangle block of Claim~\ref{s-lem:blocks}(1). Partition the \(2t\) ports
into triples and add these triples as edges. This requires \(3\mid t\).
We call the resulting hypergraphs \emph{port-matched lantern systems}.
A port triple may meet two or three blocks, with no further restriction.

\begin{proposition}\label{s-prop:equality}
For every positive integer \(t\) divisible by three, the
port-matched lantern systems constructed above are
\(\SF\)-saturated and have Berge degree four at every vertex.
They have \(5t\) vertices, \(23t/3\) edges, and deficit \(7t/3\).
They are exactly the nonempty reduced cores satisfying
\(7|V(F)|=15\deficit(F)\) in Theorem~\ref{s-thm:core}.
\end{proposition}

\begin{proof}
Let \(G\) be a port-matched lantern system with \(t\) blocks.
Each block consists of three high vertices \(a,b,c\) and
two ports \(x,y\), with edges
\(abc,abx,aby,acx,acy,bcx,bcy\).
Each high vertex has degree five, while each port has
three incident block edges and one additional port edge,
and hence has degree four. Thus these vertex classes
agree with the degree classes \(A\) and \(Q\) used above.

All links in the following argument are taken in the
whole 3-graph \(G\). For a high vertex \(a\), its link
has vertex set \(\{b,c,x,y\}\) and edge set
\(\{bc,bx,by,cx,cy\}\). Thus \(\lk_G(a)=K_4-xy\);
the missing pair consists of the two ports of its block.

For a port \(x\), let \(xuv\) be its unique additional
port edge. Its three block edges give the triangle on
\(\{a,b,c\}\) in \(\lk_G(x)\), and its additional edge
gives the link edge \(uv\). Since \(u,v\) are ports,
they are distinct from \(a,b,c\). Consequently,
\(\lk_G(x)\cong K_3\sq K_2\).
Both types of link have incidence-matching rank four,
so every vertex has Berge degree four and \(G\) is
\(\SF\)-free.

By Lemma~\ref{s-lem:sensitive}(2), every port is
universally sensitive. Hence adding any missing triple
containing a port creates a Berge five-star centred
at that port. A missing triple consisting entirely
of high vertices must meet at least two blocks,
since the three high vertices of each block already
form an edge. Adding such a triple introduces a new
neighbour at a high vertex and therefore creates a
Berge five-star, again by
Lemma~\ref{s-lem:sensitive}(2).
Thus \(G\) is saturated.

There are \(5t\) vertices, seven edges per block, and
\(2t/3\) additional port edges. Hence
\(|E(G)|=7t+2t/3=23t/3\) and
\(\deficit(G)=10t-23t/3=7t/3\).
In particular, \(7|V(G)|=15\deficit(G)\).
Since every vertex has degree at most five,
\(G\) has no \(\KFive\) component and is its own
reduced core.

Conversely, let \(F\) be a nonempty reduced core of an
\(\SF\)-saturated 3-graph, satisfying \(7|V(F)|=15\deficit(F)\).
Retain \(R,A,Q,r,a,q\) from the proof of
Theorem~\ref{s-thm:core}, and put
\(\sigma=\sum_{v\in Q}(4-d_F(v))\).
The degree count gives
\(3\deficit(F)=8r+a+2q+\sigma\), while
\(|V(F)|=5r+a+q\). Therefore
\(15\deficit(F)-7|V(F)|
=5r+(3q-2a)+5\sigma\).
All three terms on the right are nonnegative.
Equality consequently forces \(r=0\), \(2a=3q\),
and \(\sigma=0\). Thus every low vertex has degree four.

Let \(M\) be the number of shadow edges between \(A\)
and \(Q\). The adjacency estimates give
\(2a\le M\le3q\). Since \(2a=3q\), every high vertex
has exactly two low neighbours and every low vertex
has exactly three high neighbours. Each high vertex
has four neighbours in total, because its link is
\(K_4-e\). Since \(R=\varnothing\), it therefore has
exactly two high neighbours.

By Claim~\ref{s-lem:Acomponents}, every component of
\((\sh F)[A]\) is consequently a triangle.
Claim~\ref{s-lem:blocks}(1) associates with each such
triangle a seven-edge block and two ports in \(Q\).
These port sets are pairwise disjoint, since a port
of one block is adjacent to no high vertex outside
that block. They also cover \(Q\), because every
low vertex has three high neighbours.

Every port has degree three within its block and
degree four in \(F\), so it belongs to exactly one
additional edge. All edges containing a high vertex
are already among the block edges. The additional
edges therefore contain only ports and partition
the port set into triples.

Let \(t\) be the number of triangle blocks.
Since \(F\) is nonempty and the blocks cover its
vertices, \(t>0\). There are \(2t\) ports, so their
partition into triples implies \(3\mid t\).
Thus \(F\) is one of the stated port-matched
lantern systems.
\end{proof}

The standard five-lantern of~\cite{BEEJR} is the case \(t=3\) with one triple of first ports and one triple of second ports. Thus the sharp coefficient is already witnessed by a known construction, while the upper bound and equality classification are the deductions established here.

\begin{corollary}\label{s-cor:density}
A saturated 3-graph with no \(\KFive\) component has at most \(23n/15\) edges. Equality holds exactly for the empty hypergraph and the systems above.
\end{corollary}

\section{The full saturation spectrum at every order}\label{s-sec:spectrum-proof}
\subsection{The exceptional deficits near the upper endpoint}
For \(5\mid n\), the exclusions below recover~\cite[Proposition~5.1]{BEEJR}. The core bound gives a uniform treatment of every residue class and the structural descriptions at deficits zero and two.
\begin{proposition}\label{s-thm:small}
For \(\SF\)-saturated hypergraphs, the following statements hold at every order.
\begin{enumerate}
\item Deficit zero occurs exactly for disjoint unions of \(\KFive\).
\item Deficits one and three do not occur.
\item Deficit two occurs exactly for a disjoint union of copies of \(\KFive\) and one isolated vertex.
\item Deficit four occurs precisely at orders \(n\equiv1,2,4\pmod5\), with \(n\geq6,2,4\), respectively.
\end{enumerate}
The last assertion classifies feasible orders, not all isomorphism types.
\end{proposition}

\begin{proof}[Proof of Proposition~\ref{s-thm:small}]
Delete all \(\KFive\) components and write \(F=F(H)\),
\(s=|V(F)|\), and \(d=\deficit(F)=\deficit(H)\). Deletion changes
the order by a multiple of five. We retain the degree classes \(R,A,Q\)
from the proof of Theorem~\ref{s-thm:core}.

The core bound gives \(s\le\lfloor15d/7\rfloor\).
A saturated
3-graph of order at most five must be complete, since no Berge
five-star can occur on that many vertices. The complete small orders
have the following deficits:
\[
\begin{array}{c|rrrrrr}
s&0&1&2&3&4&5\\\hline
2s-\binom{s}{3}&0&2&4&5&4&0.
\end{array}
\]

For (1), if \(d=0\), then \(s=0\), so \(H\) consists entirely
of the deleted \(\KFive\) components. Conversely, every such
disjoint union is saturated and has deficit zero by
Lemma~\ref{s-lem:add}(4).

For (2), if \(d=1\), then \(s\le2\), and the table gives
no possibility. Suppose instead that \(d=3\).
Now \(s\le6\). The table excludes \(s\le5\), so suppose \(s=6\).
The stronger bound \(7s+5|R|\le45\) forces \(R=\varnothing\).
Put \(a=|A|\), \(q=|Q|\). The degree-deficiency and adjacency
counts from Theorem~\ref{s-thm:core} give \(a+q=6,\; a+2q\le9,\; 2a\le3q\). They imply \(q\le3\) and \(q\ge\lceil12/5\rceil=3\), so
\(a=q=3\). Equality in the deficiency count makes every vertex of
\(Q\) have degree four. An isolated vertex in the high shadow
would need four distinct low neighbours, but only three are
available. By Claim~\ref{s-lem:Acomponents}, the three high vertices
therefore form a triangle. Its two ports are their only low
neighbours. The third low vertex can occur only in triples contained
in \(Q\), of which there is just one, contradicting its degree four.

For (3), if \(d=2\), then \(s\le4\), and the table shows that
\(F\) is exactly one isolated vertex. Conversely, an isolated
vertex has deficit two, and adjoining any number of \(\KFive\)
components preserves saturation and deficit by Lemma~\ref{s-lem:add}(4).

For (4), set \(d=4\). We first exclude the other two residue classes.
Here \(s\le8\). Residue zero would require \(s=0\) or \(5\),
neither of which has deficit four. Residue three would require
\(s=3\) or \(8\); the table excludes \(s=3\).
If \(s=8\), then \(56+5|R|\le60\), so again \(R=\varnothing\).
The inequalities \(a+q=8,\; a+2q\le12,\; 2a\le3q\) force \(a=q=4\), with every low vertex of degree four. There are
three possibilities for the high shadow, all ruled out by the block
description in Claim~\ref{s-lem:blocks}:
\begin{enumerate}
\item If it has a triangle, the two ports of that triangle cannot be
adjacent to the fourth high vertex. That vertex is isolated in the
high shadow and needs all four low vertices as neighbours.
\item If it has an edge component, the three low neighbours of that
component cannot be shared with another high component. The remaining
high vertices have only one low vertex available, whereas each needs
at least two by Claim~\ref{s-lem:two-low}.
\item If all four high vertices are isolated, their twenty incident
edges each contain exactly one high vertex and two low vertices.
They contribute forty low incidences, but the total low degree is
only sixteen.
\end{enumerate}
The component-size partitions of four high vertices are exactly
\(3+1\), \(2+2\), \(2+1+1\), and \(1+1+1+1\).
The first case covers \(3+1\), the second covers both partitions
containing an edge component, and the third covers four isolated
vertices. Thus the classification is exhaustive and excludes \(s=8\).

To prove existence for the remaining residue classes, proceed as follows.
For residues one, two, and four, use respectively the six-vertex sun
\(S\), two isolated vertices, and \(K_4^{(3)}\). Each is saturated
and has deficit four. Adding disjoint \(\KFive\) components gives
every larger order in the same residue by Lemma~\ref{s-lem:add}(4).
The minimum orders follow from the small-order table.

\end{proof}

Proposition~\ref{s-thm:small} determines every feasible edge count
greater than \(2n-5\). It remains to construct all integer values between
\(\lceil4n/3\rceil-3\) and \(2n-5\).

\subsection{Starting hypergraphs}\label{s-sec:constructions}
The sun \(S\) was recalled in Subsection~\ref{s-sec:sun}.
The construction labels in this subsection are local to the five-star
problem. We now specify the other starting hypergraphs used by the
interval induction. The previously established constructions are cited directly; the remaining ones are verified below. The affine plane underlying \(P\) is standard, and we record its relevant saturation property without claiming the design as a new construction.

\paragraph{The 3-graph \(J\).}
Use vertices \(\{r,a_1,a_2,a_3,a_4,x,y\}\) and all edges \(ra_ia_j\) for \(i<j\), together with \(xya_i\) for \(1\leq i\leq4\). It has \((n,m,d)=(7,10,4)\). The link of \(r\) is \(K_4\); every other link is a four-edge forest. All vertices are universally sensitive and have Berge degree four.

\paragraph{The 3-graph \(X\).}
On \(\{0,\ldots,5\}\) take \(E(X)=\{012,123,234,034,014,024,135\}\).
It has \((n,m,d)=(6,7,5)\). Vertices \(1,3\) have forest links with four edges. Vertices \(0,2,4\) have Berge degree four and no neighbour \(5\). The only triples avoiding \(1,3\) are \(024\), already present, and triples containing \(5\) and two of \(0,2,4\), which introduce a new neighbour at one of those vertices. Thus \(X\) is saturated. Only vertex \(5\) has Berge degree below four. Consequently, \(X\sq K_1^{(3)}\) is also saturated: a crossing triple contains the new isolated vertex and two old vertices, at least one of which has Berge degree four.

\paragraph{The 3-graph \(Y\).}
Take \(K_4^{(3)}\) on \(\{0,1,2,3\}\) and add \(014,235\), with new vertices \(4,5\). Then \((n,m,d)=(6,6,6)\). The four old vertices have Berge degree four. Any missing triple either contains vertex \(4\) and an old vertex outside \(\{0,1\}\), or contains vertex \(5\) and an old vertex outside \(\{2,3\}\). If it contains both new vertices, at least one is a new neighbour of its old vertex. The remaining possible triples are already edges. Thus \(Y\) is saturated.

\paragraph{The affine-plane 3-graph \(P\).}
Take the twelve lines of the affine plane of order three on nine vertices. Equivalently use
\(
012,345,678,036,147,258,048,156,237,057,138,246.
\)
Each vertex has link \(4K_2\). Thus \(P\) is saturated and attachment-ready, with \((n,m,d)=(9,12,6)\).

\paragraph{The 3-graph \(U\).}
Take the eleven edges \(701,702,703,712,713,014,245,256,346,356,456\) on \(\{0,\ldots,7\}\). It has \((n,m,d)=(8,11,5)\). Vertices \(0,1,2,3,4\) have four-edge forest links and are universally sensitive. Vertices \(5,6\) have four-edge triangle-with-pendant-edge links, and vertex \(7\) has link \(K_4-23\); all have Berge degree four. The only triple avoiding the five sensitive vertices is \(567\). Adding it gives vertex \(5\) the new neighbour \(7\). Thus \(U\) is saturated and attachment-ready.

\paragraph{The 3-graph \(V\).}
Take a triangle block on high vertices \(0,1,2\) and ports \(3,4\), a complete 3-graph on \(\{5,6,7,8\}\), and edges \(356,478\). Then \((n,m,d)=(9,13,5)\). Ports \(3,4\) have links \(K_3\sq K_2\) and are universally sensitive. All remaining vertices have Berge degree four. Every triple contained in either the high triple or the complete four-set is already an edge. Thus every missing triple avoiding both ports meets both sets and, when added, gives a high vertex a new neighbour. Thus \(V\) is saturated and attachment-ready.

\paragraph{The 3-graph \(Z\).}
Take the Fano plane on \(\{0,\ldots,6\}\) with edges
\(013,124,235,346,045,156,026\), and add \(012,345\). It has \((n,m,d)=(7,9,5)\). The six vertices in the added triples have four-edge forest links: their original links are three disjoint edges, and the new edge joins two components. They are universally sensitive. The sole remaining vertex has degree three, so every missing triple meets a sensitive vertex. Hence \(Z\) is saturated.

\paragraph{The 3-graph \(W\).}
Take \(K_4^{(3)}\) on \(\{0,1,2,3\}\), then add \(045,145,246,347,567\).
It has \((n,m,d)=(8,9,7)\). The four old vertices have links \(K_3\sq K_2\), and vertex \(4\) has a four-edge forest link. The remaining three vertices form the edge \(567\). Thus every missing triple meets a universally sensitive vertex, proving saturation.

\paragraph{The 3-graph \(T\).}
Start with the sun whose edges are \(4uv,5uv\) for \(uv\in\{67,78,89,96\}\), delete \(467,589\), and add
\(
 024,135,267,389,
 \{10,0,1\},\{10,0,2\},\{10,0,3\},\{10,1,2\},\{10,1,3\}.
\)
The vertex set is \(\{0,\ldots,10\}\) and \((n,m,d)=(11,15,7)\). Every vertex except \(10\) has a four-edge forest link. To see this directly, each vertex of the sun loses one edge from its four-cycle link and gains one edge containing a vertex of \(\{0,1,2,3\}\); the links at \(0,1,2,3\) are also four-edge forests. Vertex \(10\) has link \(K_4-23\). Therefore every missing triple meets a sensitive vertex. The 3-graph is saturated and attachment-ready.

\paragraph{The known 3-graphs \(B\) and \(D\).}
Let \(B\) and \(D\) be the broken lantern and the 3-graph in
\cite[Constructions~5.8 and~5.9]{BEEJR}, respectively.
Their parameters \((n,m,d)\) are \((10,15,5)\) and \((10,14,6)\).
Both are saturated by~\cite[Proposition~5.12]{BEEJR}; these
are the only properties of \(B\) and \(D\) needed below.

\subsection{Starting orders and completion of the spectrum}
Table~\ref{s-tab:initial} supplies every integer in \([\lceil4n/3\rceil-3,2n-5]\) for \(7\leq n\leq12\). In each row, the constructions are listed in the same order as the consecutive edge counts. Each union has an attachment-ready saturated factor; the sole exception \(X\sq C_1\) has already been proved saturated.

\begin{table}[ht]
\centering
\small
\renewcommand{\arraystretch}{1.45}
\begin{tabular}{ccl}
\toprule
\(n\)&edge counts&constructions, in increasing edge-count order\\\midrule
7&\(7,8,9\)&\(X\sq C_1,\ S\sq C_1,\ Z\)\\
8&\(8,9,10,11\)&\(S\sq C_2,\ W,\ J\sq C_1,\ U\)\\
9&\(9,10,11,12,13\)&\(S\sq C_3,\ J\sq C_2,\ U\sq C_1,\ P,\ V\)\\
10&\(11,12,13,14,15\)&\(J\sq C_3,\ P\sq C_1,\ V\sq C_1,\ D,\ B\)\\
11&\(12,13,14,15,16,17\)&\(P\sq C_2,\ V\sq C_2,\ J\sq C_4,\ T,\ Y\sq C_5,\ X\sq C_5\)\\
12&\(13,14,15,16,17,18,19\)&\(P\sq C_3,\ V\sq C_3,\ U\sq C_4,\ 2S\),\\
&&\(X\sq C_1\sq C_5,\ S\sq C_1\sq C_5,\ Z\sq C_5\)\\\bottomrule
\end{tabular}
\caption{Starting 3-graphs for an induction on the order.}\label{s-tab:initial}
\end{table}

\begin{proof}[Proof of Theorem~\ref{s-thm:spectrum}]
We first show by induction that every integer in
\([\lceil4n/3\rceil-3,2n-5]\) occurs for every \(n\geq7\).
The starting orders \(7\leq n\leq12\) are covered by
Table~\ref{s-tab:initial}. Suppose \(n\geq13\) and apply induction
at orders \(n-6\) and \(n-5\). Adding the sun \(S\), which has six
vertices and eight edges and is attachment-ready, gives every integer in
\[
 [\lceil4(n-6)/3\rceil-3+8,\,2(n-6)-5+8]
 =[\lceil4n/3\rceil-3,\,2n-9].
\]
Adding \(\KFive\) gives every integer in
\[
 [\lceil4(n-5)/3\rceil-3+10,\,2(n-5)-5+10]
 =[\lceil(4n+1)/3\rceil,2n-5].
\]
For \(n\geq13\), we have \(\lceil(4n+1)/3\rceil\leq2n-8\),
since \(4n+1\leq6n-24\). The two integer intervals overlap or are
adjacent and together fill \([\lceil4n/3\rceil-3,2n-5]\), completing
the induction.

It remains to identify the upper endpoint exceptions and exclude all
other values. Proposition~\ref{s-thm:small} gives all possibilities
above \(2n-5\): deficit zero only in residue zero; deficit two only in
residue one; and deficit four only in residues one, two, and four.
Deficits one and three are absent. Lemma~\ref{s-lem:sat-lower}(5)
excludes all lower edge counts, and Lemma~\ref{s-lem:local}(3)
excludes edge counts above \(2n\).

At order six, the 3-graphs \(Y,X,S\) give respectively six, seven, and eight edges, and \(\KFive\sq C_1\) gives ten. Lemma~\ref{s-lem:sat-lower}(5) excludes smaller values, while Proposition~\ref{s-thm:small} excludes nine and eleven by deficits three and one, and excludes twelve because deficit zero requires an order divisible by five. Lemma~\ref{s-lem:local}(3) excludes more than twelve edges. At orders at most five, no five-star can occur even after adding a triple, so saturation requires the complete 3-graph. This proves the theorem at every order.
\end{proof}

\begin{proof}[Proof of Corollaries~\ref{s-thm:saturation} and~\ref{s-thm:extremal}]
Take the minimum and maximum of the spectra in Theorem~\ref{s-thm:spectrum}. The respective upper endpoints in residues \(0,1,2,3,4\) modulo five are \(2n,2n-2,2n-4,2n-5,2n-4\), which equal \(10q+\binom t3\) for \(n=5q+t\). The stated construction \(q\KFive\sq C_t\) is free and has that many edges. The small orders follow from the same theorem.
\end{proof}

\begin{remark}
The proof gives an explicit construction for every feasible pair \((n,m)\).
For a central-interval value with \(n\ge13\), subtract \((6,8)\) when
\(m\le2n-9\), and subtract \((5,10)\) otherwise. The interval calculation
keeps the residual pair in the required range until a row of
Table~\ref{s-tab:initial} applies. Restoring the deleted sun and complete
components then produces the desired hypergraph.
\end{remark}

\section{Local tools for Berge-\texorpdfstring{\(K_4\)}{K4} saturation}\label{k-sec:local}
Throughout this section and the remainder of the paper, saturation and
freeness refer to \(\KF\). We use the local framework of~\cite{CHX}
and include the local arguments needed for the improved lower bound.

In a \(\KF\)-free 3-graph, a pair \(uv\) is \emph{good} if there are
distinct vertices \(s,t\notin\{u,v\}\) and five distinct hyperedges
\(e_{us},e_{ut},e_{vs},e_{vt},e_{st}\) representing
\(us,ut,vs,vt,st\), respectively. Otherwise \(uv\) is \emph{bad}.
This ordered representation, together with \(u,v,s,t\), is called a
\emph{witness of Berge \(K_4\)} for \(uv\). It represents
\(K_4-uv\); a new hyperedge containing \(uv\) supplies the sixth
representative and completes a Berge \(K_4\).
The vertices \(s,t\) are its \emph{additional core vertices}.
Each has degree at least three, since three distinct representative
hyperedges contain it. The representative \(e_{st}\) of the pair
joining the additional core vertices is called the \emph{closure}
of the witness of Berge \(K_4\).

A \(\KF\)-free 3-graph is \(\KF\)-saturated if and only if
every missing triple contains a good pair. Indeed, a witness of
Berge \(K_4\) for a pair in that triple, together with the new
hyperedge, gives the forbidden copy. Conversely, every forbidden
copy after an addition uses the new hyperedge. Deleting its
representative leaves a witness of Berge \(K_4\) for the pair
represented by that hyperedge.
Whenever a witness of Berge \(K_4\) is displayed as an ordered list,
the representatives are listed in the order
\(us,ut,vs,vt,st\) for the indicated pair \(uv\) and additional core
vertices \((s,t)\). In contrast, a witness of Berge five-star,
as defined in Subsection~\ref{s-sec:preliminaries}, completes four old
representatives to a five-edge star centred at a specified vertex.

\begin{lemma}\label{k-lem:goodfacts}\label{k-lem:same-centre}\label{k-lem:mindegree}
The following hold in every \(\KF\)-saturated 3-graph of order at least five.
\begin{enumerate}
\item If \(uv\) is good and \(uvw\in E(H)\), every witness of
Berge \(K_4\) for \(uv\) uses \(uvw\), and \(w\) is an additional core vertex.
Thus \(d(w)\ge3\), and a good pair has codegree at most two.
\item All pairs in an edge containing two vertices of degree
at most two are bad.
\item A vertex in an edge whose three vertices have degree
at most two has no good partner.
\item If \(d(p)=2\) and \(pv\) is good, the two additional
core vertices in any witness of Berge \(K_4\) for \(pv\) are
distinct neighbours of \(p\) of degree at least three. The two edges at \(p\) represent its two old
core pairs.
\item There are no degree-two vertices \(p,q,r\), with \(q\ne r\),
such that the edges at \(p\) are \(pqz,prz\).
\item If \(n\ge7\) and \(H\) has the minimum number of edges among
all \(n\)-vertex \(\KF\)-saturated 3-graphs, then \(\delta(H)\ge2\).
\end{enumerate}
\end{lemma}
\begin{proof}
For (1), let \(\mathcal W\) be a witness of Berge \(K_4\) for
\(uv\). If \(uvw\) were not one of its five representatives,
it could represent \(uv\), giving a forbidden Berge \(K_4\).
Thus \(uvw\) is used. Each of the five pairs represented by \(\mathcal W\)
involves an additional core vertex, so \(w\)
is one of the two additional core vertices. Each additional core
vertex belongs to three distinct representatives, giving \(d(w)\ge3\).
The same argument applies to every edge containing \(uv\), and
there are only two possible additional core vertices; hence
\(d(u,v)\le2\).

For (2), let \(pqr\) be an edge with \(d(p),d(q)\le2\).
By (1), the pairs \(pr,qr\) are bad. If \(pq\) were good,
its witness of Berge \(K_4\) would use \(pqr\) for a pair
incident with exactly one of \(p,q\). The other endpoint would
then need two representatives different from \(pqr\), although
it has at most one other incident edge. This is impossible.

For (3), let \(pqr\) be an edge with all three degrees at most two.
The pairs \(pq,pr\) are bad by (2). A witness of Berge \(K_4\)
for \(pv\), where \(v\notin\{p,q,r\}\), would have to use
\(pqr\) as one of the two representatives at \(p\). It would
therefore make \(q\) or \(r\) an additional core vertex of degree
at least three, a contradiction. These three assertions are also
Fact~3.1 of~\cite{CHX}.

For (4), the two pairs incident with \(p\) in a witness of Berge
\(K_4\) require distinct representatives. Since \(d(p)=2\), these
are exactly its two incident edges. The additional core vertices
are consequently distinct neighbours of \(p\), and their degrees
are at least three by the definition of the witness of Berge \(K_4\).

For (5), suppose the two edges at \(p\) are \(pqz,prz\).
The vertex \(p\) has at most one neighbour of degree at least three,
so the fourth assertion gives it no good partner. Saturation of the
missing triples through \(p\) then makes every pair outside \(p\),
except possibly \(qz,rz\), good. The other edge at \(q\) avoids \(p\),
and its pair opposite \(q\) is bad by the first assertion. It must
therefore be \(qrz\). But \(qr\) is bad by the second assertion,
contradicting the preceding description of bad pairs outside \(p\).
For (6), apply Lemma~3.3 of~\cite{CHX}, whose minimum-edge
assumption is exactly the assumption in this assertion.
\end{proof}

For a \(\KF\)-saturated 3-graph \(H\) with \(\delta(H)\ge2\), put \(\Low=\{v:d(v)=2\}\) and
\(\High=\{v:d(v)\ge3\}\), and write \(\lambda=|\Low|\), \(x=|\High|\),
and \(q_{\High}=|E(H[\High])|\). Vertices in \(\Low\) and \(\High\)
are called \emph{low} and
\emph{high}, respectively, and \(N_{\High}(v)=N_H(v)\cap\High\).
These degree classes apply only to the Berge-\(K_4\) part and differ
from the classes \(A,Q\) in the five-star part.
If \(\Low\ne\varnothing\), then \(n\ge5\): a saturated
3-graph of order at most four is complete and has no degree-two vertex.
No edge is contained in \(\Low\). Indeed, if \(pqr\) were such an
edge, each of \(p,q,r\) would have no good partner by
Lemma~\ref{k-lem:goodfacts}(3). Choose \(z\notin\{p,q,r\}\).
Saturation would force both \(pqz\) and \(prz\) to be edges,
contradicting \(d_H(p)=2\).
An edge meeting \(\Low\) in two vertices has \emph{type~1};
one meeting it in one vertex has \emph{type~2}.
The \(\High\)-pair of a type~2 edge is its \emph{support pair}
(or \emph{edge support}). Edges contained in \(\High\) are called
\emph{all-high edges}. The \emph{vertex support} of a low vertex
\(v\) is \(N_{\High}(v)\).

Let \(H\) be an \(n\)-vertex, \(m\)-edge
\(\KF\)-saturated 3-graph with \(\delta(H)\ge2\),
and let \(p,q\in\Low\) form a bad pair.
Set \(S=N_{\High}(p)\cup N_{\High}(q)\), \(s=|S|\),
\(U=V(H)\setminus(\{p,q\}\cup S)\), and
\(I=\{e\in E(H):e\cap\{p,q\}\ne\varnothing\}\).
Write \(D(S)=\sum_{v\in S}d_H(v)\) and
\(\kappa=\sum_{e\in I}(|e\cap S|-\tfrac12|e\cap U|)\).

For \(i=2,3\), put
\(h_i=|\{e\in E(H)\setminus I:|e\cap S|=i\}|\).
Project each edge in
\(\mathcal P=\{e\in E(H)\setminus I:|e\cap S|\ge2\}\)
onto \(S\), retaining its edge identity.
Give its projection weight \(w(e)=|e\cap S|-1\).
For \(\mathcal R\subseteq\mathcal P\), write
\(w(\mathcal R)=\sum_{e\in\mathcal R}w(e)\).
Thus the total projection weight is
\(W=w(\mathcal P)=h_2+2h_3\).

\begin{lemma}\label{k-lem:charge}\label{k-lem:twotriangles}
With the notation above, the following hold.
\begin{enumerate}
\item
\(D(S)\ge n-s-2+\kappa+\tfrac32W\), and
\(3m\ge2\lambda+3(x-s)+D(S)\).

\item
If \(\mathcal R_i\subseteq\mathcal P\) is the set of the three
representatives of a Berge triangle with core \(C_i\subseteq S\), then
\(w(\mathcal R_i)\ge3\).
For \(C_1\ne C_2\),
\[
w(\mathcal R_1\cup\mathcal R_2)\ge
\begin{cases}
5,& |C_1\cap C_2|=2,\\
6,& |C_1\cap C_2|\le1.
\end{cases}
\]
\end{enumerate}
\end{lemma}

\begin{proof}
For (1), the triple \(pqv\) is missing for every \(v\in U\):
otherwise \(v\) would be low by the definition of \(S\), giving
an edge contained in \(\Low\).

Since \(pq\) is bad, saturation of \(pqv\) forces
\(pv\) or \(qv\) to be good. Choose a witness of
Berge \(K_4\) for such a pair.
By Lemma~\ref{k-lem:goodfacts}(4), its two additional
core vertices belong to \(S\).
Mark the two distinct representatives joining \(v\)
to these vertices. Let \(\mu(e)\) count the marks
received by \(e\). Then
\(\sum_{e\in E(H)}\mu(e)=2|U|\).
Each vertex of \(U\) marks a given edge at most once,
so \(\mu(e)\le|e\cap U|\).
Moreover, an edge disjoint from \(S\) receives no marks.
Call \(|e\cap S|-\mu(e)/2\) the \emph{surplus} of \(e\)
in this count.

For \(e\notin I\) with \(k=|e\cap S|\ge1\), we have
\(\mu(e)\le3-k\), and hence
\(k-\tfrac12\mu(e)\ge\tfrac32(k-1)\).
For \(e\in I\), we have
\(|e\cap S|-\tfrac12\mu(e)
\ge |e\cap S|-\tfrac12|e\cap U|\).
Summing over all edges yields
\(D(S)-|U|\ge\kappa+\tfrac32h_2+3h_3
=\kappa+\tfrac32W\).
Since \(|U|=n-s-2\), this proves the first inequality.
The \(\lambda\) low vertices have degree two,
and the \(x-s\) high vertices outside \(S\)
have degree at least three. Summing their degrees
together with \(D(S)\) proves the second inequality.

For (2), a Berge triangle uses three distinct
representatives, each of weight at least one.
Thus \(w(\mathcal R_i)\ge3\).

Put \(c=|\mathcal R_1\cap\mathcal R_2|\).
The union contains \(6-c\) distinct representatives.
A shared projection of weight one must represent
a pair common to both cores.

If \(|C_1\cap C_2|=2\), there is only one common
core pair, and each triangle uses only one
representative for it. Hence at most one shared
projection has weight one.
For \(c\ge1\), at least \(c-1\) shared projections
therefore have weight two, giving
\(w(\mathcal R_1\cup\mathcal R_2)
\ge(6-c)+(c-1)=5\).
For \(c=0\), the union weight is at least six.

If \(|C_1\cap C_2|\le1\), there is no common
core pair, so every shared projection has weight two.
Consequently,
\(w(\mathcal R_1\cup\mathcal R_2)\ge(6-c)+c=6\).
\end{proof}

\section{The lower bound for Berge-\texorpdfstring{\(K_4\)}{K4}}\label{k-sec:lower}
\subsection{Edges containing two degree-two vertices}
We adapt the low-vertex reductions of~\cite[Section~3.2]{CHX}
to the full order range. The short boundary arguments for the known
cases \(n=7,8\) are retained so that the lower-bound proof is analytic
at every order. Throughout this section, a \emph{counterexample} means an
\(n\)-vertex \(\KF\)-saturated 3-graph \(H\) with \(n\ge7\) and
\(m=\sat_3(n,\KF)\le n-1\). This minimum-edge assumption makes
Lemma~\ref{k-lem:goodfacts}(6) applicable. Hence the partition
\(\Low\cup\High\) is defined, and degree summation gives
\(3n-3\ge2\lambda+3x=3n-\lambda\), so \(\lambda\ge3\).
Only the stated degree and saturation properties are used below.

\begin{proposition}\label{k-prop:typeone}\label{k-lem:forbidden}\label{k-lem:decompose}\label{k-lem:matching-charge}
In a counterexample, the low pairs belonging to type~1 edges form a
matching. If a type~1 edge exists, then \(x\le5\). Moreover,
\(n\le3x-2\), with the stronger bound \(n\le2x-1\) when there
are no type~1 edges.
\end{proposition}
\begin{proof}
Form a loopless multigraph \(M\) on \(\Low\), placing one edge
\(pq\) for each type~1 hyperedge \(pqh\), where \(h\in\High\).
At this stage distinct hyperedges may give parallel edges of \(M\).
Since every low vertex has hypergraph degree two, every vertex of
\(M\) has degree at most two. We first restrict its components,
then exclude the only possible components with two edges.

\begin{claim}
The multigraph \(M\) is simple, and each of its components is an
isolated vertex, an isolated edge, or a path with two edges.
Consequently, under the counterexample assumption \(m\le n-1\),
we have \(n\le4x-3\).
\end{claim}
\begin{claimproof}
Suppose first that two low vertices \(p,q\) have incident edges
\(pqz,pqw\). The missing triple \(pzw\) forces \(zw\) good, so
\(x\ge4\). For the bad pair \(pq\), Lemma~\ref{k-lem:charge} has
\(S=\{z,w\}\) and \(\kappa=2\), initially giving \(D(S)\ge n-2\).
Every witness of Berge \(K_4\) forced by \(pqv\), \(v\in U\),
has additional core vertices
\(z,w\); hence some edge outside \(I\) contains \(zw\). Two such
edges contribute surplus at least three in the incidence count,
giving \(D(S)\ge n+1\). If there is only one, say \(zwv\), then \(v\in U\).
In the witness of Berge \(K_4\) forced by the missing triple
\(pqv\), this edge must represent \(zw\). It therefore cannot
also represent either marked pair at \(v\). Since \(v\) is its
only vertex in \(U\), the edge receives no marks. Its surplus is two, so \(D(S)\ge n\). The degree bound in
Lemma~\ref{k-lem:charge} now gives \(3m\ge3n+x-6\), hence
\(x\le3\), a contradiction. Thus \(M\) has no parallel edges.

If \(p\) has edges \(pqz,prw\), Lemma~\ref{k-lem:goodfacts}(5)
gives \(z\ne w\). The missing triple \(pqr\) forces \(qr\) good
with additional core vertices \(z,w\). The other edges at \(q,r\)
therefore contain \(w,z\), respectively. The missing triple
\(pqw\) forces \(qw\) good, so its other edge is \(qwt\) with
\(t\in\High\setminus\{z,w\}\). Thus \(q\), and similarly \(r\),
has only one type~1 edge, proving the component description.
To obtain the numerical consequence, let \(i,j,k\) count the
components of \(M\) with one, two, and three vertices, respectively.
An isolated vertex of \(M\) still has degree two in \(H\), so it
belongs to two type~2 edges. An isolated edge of \(M\) accounts for
one type~1 edge and two type~2 edges. A two-edge path accounts for
two type~1 edges and two type~2 edges. No hyperedge is counted for two distinct components: an edge meeting low vertices
in two components would join those components in \(M\).
Consequently, \(\lambda=i+2j+3k, \; m=2i+3j+4k+q_{\High}\ge\frac43\lambda\). Combining this with \(m\le\lambda+x-1\) gives \(\lambda\le3x-3\), and hence
\(n=\lambda+x\le4x-3\). The same identity gives \(m\ge3\lambda/2\) when
\(k=0\), and \(m\ge2\lambda\) when \(j=k=0\).
\end{claimproof}

We next exclude the two-edge paths. For a path centre \(p\) and
endpoint \(q\), their incident edges have the form \(pqz,prw,qwt\),
with \(S=\{z,w,t\}\) and \(\kappa=7/2\). The good pair \(qw\)
forced by \(pqw\) gives three triangle representatives outside \(I\),
so Lemma~\ref{k-lem:charge} yields \(D(S)\ge n+3\). Its degree
bound gives \(x\le3\), hence \(\High=S\). The triple \(zwt\) must be an edge: adding it leaves every low
vertex at degree two, whereas all four core vertices of a Berge
\(K_4\) must have degree at least three. Only the three vertices
of \(\High=S\) could be core vertices. The preceding count, including this all-high
edge, gives \(m\ge4\lambda/3+1\), whereas \(m\le\lambda+2\).
Thus \(\lambda\le3\) and \(n\le6\), a contradiction.
The type~1 pairs therefore form a matching. The two sharper counts
in the claim give \(n\le3x-2\) and, in the absence of type~1 edges,
\(n\le2x-1\).

It remains to bound \(x\) when a matching pair exists. We record the
required incidence estimate within the proof.

\begin{claim}
Suppose the edges at two low vertices are \(pqz,puv,qwt\), where
\(z,u,v,w,t\in\High\). For \(S=N_{\High}(p)\cup N_{\High}(q)\),
\[
 D(S)\ge
 \begin{cases}
 n+6-s,&s\le3,\\
 n+4,&s=4,\\
 n+7,&s=5.
 \end{cases}
\]
\end{claim}
\begin{claimproof}
Here \(\kappa=5\). For every \(r\in S\setminus\{z\}\), the
missing triple \(pqr\) forces a good pair and a Berge triangle
containing \(z,r\) in \(S\). Its representatives avoid the two
used edges at the selected low endpoint; only the other endpoint's
type~2 edge can be borrowed. Thus \(W\ge2\), proving the first bound
by Lemma~\ref{k-lem:charge}.

For \(s=4\), if \(z\) lies in neither type~2 support, write the
supports as \(uv,uw\). The triangle forced by \(pqu\) is represented
entirely outside \(I\). If \(z\) belongs to a support, write the
edges as \(pqz,pzu,qwt\); the same holds for the triangle forced
by \(pqw\). Hence \(W\ge3\), and integrality gives \(D(S)\ge n+4\).

For \(s=5\), the four vertices other than \(z\) are distinct and every
selected triangle is represented outside \(I\). The pairs opposite
\(z\) in these triangles cover those four vertices. A matching of size two gives \(W\ge6\)
by the triangle-weight bound in Lemma~\ref{k-lem:charge}. Otherwise
the pairs form a three-edge star, so three triangle cores share one
pair and have three distinct remaining vertices. Their pair union
has seven pairs. Weight at most five with zero or one triple
projection covers at most five or six pairs. With two triple
projections there is at most one pair projection, and each triangle
must use all three edges. The pair projection must be the common
pair; each triple projection would then have to contain the three
remaining vertices and an endpoint of that pair, which is impossible.
Thus \(W\ge6\), giving the last bound.
\end{claimproof}

Apply this claim to a type~1 pair and use the degree bound in
Lemma~\ref{k-lem:charge}. For \(s\le3\), it gives \(x\le3\),
while for \(s=4,5\) it gives \(x\le5\). This proves the proposition.
\end{proof}

\begin{lemma}\label{k-lem:twohigh}\label{k-lem:two-high-not-disjoint}
For any saturated 3-graph \(H\) with \(\delta(H)\ge2\), the following statements hold.
\begin{enumerate}
\item Suppose \(n\ge7\), \(m\le n\), and the type~1 pairs form a
matching. Then \(q_{\High}\ge1\). If \(q_{\High}=1\), necessarily
\(n=m=7\), \(x=3\), \(\lambda=4\), and the matching has two edges.
In particular, \(q_{\High}\ge2\) if \(m\le n-1\) or if there
are no type~1 edges.
\item If there are exactly two all-high edges, they intersect.
\end{enumerate}
\end{lemma}
\begin{proof}
For (1), the matching count gives
\(m\ge3\lambda/2+q_{\High}\), so \(n\le3x\)
when \(q_{\High}=0\); in particular \(x\ge3\).
Every positive-codegree pair in \(\High\) is then bad, since its
representatives have low third vertices. A good pair in \(\High\)
would yield an old Berge triangle whose core is a missing
triple of bad pairs. Hence all pairs in \(\High\) are bad,
contradicting saturation of a triple contained in \(\High\).

Suppose \(q_{\High}=1\), with edge \(P=abc\). A positive-codegree
good pair in \(\High\) must lie in \(P\), which is the only
hyperedge containing that pair. If \(ab\) were good, its additional core
vertices would be \(c,d\), with \(d\notin P\). Its old
triangles on \(acd,bcd\) have bad type~2 pairs at \(d\);
saturation forces \(ac,bc\) to be good. Then \(P\) would be the
only hyperedge containing either pair, although they require distinct
representatives in the witness of Berge \(K_4\) for \(ab\). Thus all positive-codegree pairs are bad.
A zero-codegree good pair yields two distinct old triangle
cores of bad pairs, both of which would have to equal \(P\).
So all pairs in \(\High\) are bad and \(x\le3\).
Now \(3\lambda/2+1\le m\le\lambda+x\) gives \(\lambda\le2x-2\).
The assumption \(n\ge7\) forces \(x=3,\lambda=4,m=7\), and
the identity \(m=2\lambda-k+q_{\High}\), where \(k\) is the matching size,
gives \(k=2\). This also proves the stated consequences.

For (2), suppose the two all-high edges \(P,Q\)
are disjoint. A positive-codegree good pair
\(ab\subseteq P=abc\) is contained in no hyperedge other than
\(P\), and any witness of Berge \(K_4\) for \(ab\) has core
\(\{a,b,c,d\}\), where \(d\notin P\). In its two old triangles
\(acd,bcd\), the pairs at \(d\) have type~2 representatives
and are bad. Saturation forces \(ac,bc\) to be good, so \(P\)
would be the only hyperedge containing either pair. These two pairs
require distinct representatives in the witness of Berge \(K_4\)
for \(ab\), a contradiction. Thus all
positive-codegree pairs in \(\High\) are bad. A zero-codegree
good pair gives two distinct old triangle cores of bad pairs;
they must be \(P,Q\), yet share two vertices. Hence all pairs
in \(\High\) are bad, contradicting a missing triple contained in \(\High\).
\end{proof}

\begin{proposition}\label{k-prop:typeone-orders-eight-nine}\label{k-cor:no-type-one}\label{k-cor:typeone-nine}\label{k-lem:five-high-two-triples}
There is no counterexample of order eight. Every counterexample of
order \(n\ge8\) has no type~1 edge.
\end{proposition}
\begin{proof}
By Proposition~\ref{k-prop:typeone}, the type~1 pairs form a matching.
If its size is \(k\), Lemma~\ref{k-lem:twohigh} gives
\(m=2\lambda-k+q_{\High}\ge3\lambda/2+2\). Hence \(\lambda\le2x-6\) and
\(n\le3x-6\). If a type~1 edge exists, \(x\le5\), so \(n\le9\).
We exclude the remaining orders eight and nine using the following
local observation.

\begin{claim}
A counterexample cannot have \(x=5\), two all-high edges, and at
most four type~2 edges.
\end{claim}
\begin{claimproof}
Let the all-high edges be \(P,Q\). A witness of Berge \(K_4\) for a pair in \(\High\) uses only
all-high and type~2 edges. Positive-codegree good pairs
lie in \(P\) or \(Q\) and are not type~2 supports.

First let \(P=abc,Q=ade\). If \(bc\) were good, its core,
after interchanging \(d,e\) if necessary, would be \(\{a,b,c,d\}\), forcing type~2 supports
\(bd,cd\) and at least one of \(ab,ac\), say \(ab\).
Missing \(abd\) then forces \(ad\) good. Its core is
\(\{a,d,e,t\}\), \(t\in\{b,c\}\), forcing supports \(dt,et\)
and one of \(ae,de\). These two collections share only
\(dt\), giving five supports, a contradiction. Thus
\(bc,de\) are bad. A good pair \(ab\) would give an old
triangle on \(bct\), \(t\in\{d,e\}\), whose three pairs
are bad and whose triple is missing. This also excludes
\(ac,ad,ae\). All positive-codegree pairs are therefore bad.
A zero-codegree good pair gives two old triangle cores which
must be \(P,Q\), but they share two vertices, whereas
\(|P\cap Q|=1\). Hence every pair in \(\High\) is bad,
contradicting saturation.

Now let \(P=abc,Q=abd\), and \(\High=\{a,b,c,d,e\}\).
No witness of Berge \(K_4\) for a pair in \(\High\) has core
\(\{a,b,c,d\}\): it would use at least
three type~2 edges supported on these vertices, leaving at most one
type~2 edge incident with \(e\). Thus \(e\) would have no
good partner in \(\High\), and a supported pair \(uv\subseteq\{a,b,c,d\}\)
would give the missing triple \(euv\) of bad pairs.
In particular \(ab\) is bad, since \(P,Q\) force its core
to be \(\{a,b,c,d\}\). A good pair contained in exactly one of
\(P,Q\), say \(ac\), must have core \(\{a,b,c,e\}\); its triangle \(abe\) is missing and has
three bad pairs. The same excludes \(bc,ad,bd\).
Finally, a zero-codegree good pair would have its two old
triangle cores equal to \(P,Q\), giving the excluded core
\(\{a,b,c,d\}\). Again all pairs in \(\High\) are bad.
\end{claimproof}

We have \(\lambda\ge3\), \(q_{\High}\ge2\), \(k\le\lfloor\lambda/2\rfloor\),
and \(m=2\lambda-k+q_{\High}\). For \(n=8\), these force
\((\lambda,x,m,k,q_{\High})=(3,5,7,1,2)\). For \(n=9\) with a type~1
edge, \(x\le5\) forces \((\lambda,x,m,k,q_{\High})=(4,5,8,2,2)\).
Both cases have four type~2 edges and contradict the claim.
There is consequently no counterexample of order eight, and the
bound \(n\le9\) excludes type~1 edges at every order \(n\ge8\).
\end{proof}

\subsection{Bad pairs of degree-two vertices}\label{k-sec:bad-low}
Suppose there are no type~1 edges and \(p,q\in\Low\) form a
bad pair. Use the notation \(S,s,I,h_2,h_3\) introduced before
Lemma~\ref{k-lem:charge}.
The four edges of \(I\) are distinct type~2 edges, so \(\kappa=8\)
and \(s\le8\). Each of \(N_{\High}(p)\) and \(N_{\High}(q)\)
has three or four vertices. Call the vertices of \(S\)
\emph{support vertices}. Project onto \(S\) every edge outside \(I\)
meeting \(S\) at least twice, retaining its edge identity, to obtain
a hypergraph \(J\) on \(S\) with edges of sizes two and three.
Equal pair projections from different original edges remain distinct.
The degree \(d_J(v)\) counts incident edge identities, and the total
weight is \(W=h_2+2h_3\). An \emph{incidence component} of \(J\)
means a connected component of its incidence graph.

A Berge triangle with core in \(S\) is called \emph{complete}
if all three of its representatives avoid \(I\). Such triangles
correspond to Berge triangles in \(J\), and hence to six-cycles
in its simple bipartite incidence graph.
For each \(v\in S\), the missing triple \(pqv\) forces \(pv\)
or \(qv\) good. A witness of Berge \(K_4\) for that pair contains
a Berge triangle on \(v\) and its two additional core vertices,
all in \(S\). We call any such represented triangle a
\emph{triangle forced at \(v\)}.

Combining the incidence inequality with the degree sum gives
\(3m\ge2\lambda+3(x-s)+n-s+6+\lceil3W/2\rceil\). Since \(\lambda=n-x\),
this rearranges to \(x\le4s-6+3(m-n)-\lceil3W/2\rceil\).
In a counterexample, therefore,
\(x\le4s-9-\lceil3W/2\rceil\).
The absence of type~1 edges also gives \(m=2\lambda+q_{\High}\), whence
\(\lambda\le x-q_{\High}-1\) and \(n\le2x-3\). Here \(\lambda\ge3\) and
\(q_{\High}\ge2\) by Lemma~\ref{k-lem:twohigh}, so \(x\ge6\).

The following cycle-space estimate supplies the only graph-theoretic
input needed for the support analysis.

\begin{theorem}\label{k-thm:short-cycle-support}\label{k-lem:cyclecover}
Let \(B\) be a finite simple bipartite graph with a specified
vertex class \(V_0\), and let \(\ell\ge2\) be an integer. Suppose a family
\(\mathcal C\) of distinct simple cycles of length at most
\(2\ell\) spans a subspace of dimension \(d\) over
\(\mathbb F_2\). If \(|\mathcal C|\ge d+1\), then
\(
 \left|V_0\cap\bigcup_{C\in\mathcal C}V(C)\right|
 \le\ell d-\left\lceil\frac\ell2\right\rceil.
\)
If the family spans the full cycle space, the bound covers
every vertex of \(V_0\) lying on a cycle. It is sharp for
every \(\ell\ge2,d\ge2\), even among connected graphs
whose simple cycles all have length \(2\ell\).

In particular, a connected bipartite graph of cycle rank two has at
most three simple cycles; if it has three six-cycles, at most four
vertices in either class lie on cycles. At cycle rank three, four
distinct six-cycles give the corresponding bound seven.
\end{theorem}
\begin{proof}
Choose \(C_1,\ldots,C_d\) from \(\mathcal C\) as a basis of
its span, and choose one additional cycle from \(\mathcal C\).
The dependence expressing this additional cycle in the basis involves
\(k\ge3\) cycles. Each
vertex in that dependence occurs in at least two cycles;
otherwise an incident edge occurs only once in the zero sum.
The dependence therefore covers at most
\(\lfloor k\ell/2\rfloor\) vertices of \(V_0\), and the
remaining selected cycles contribute at most \(\ell\) each.
The full union has at most
\[
 \ell(d+1-k)+\lfloor k\ell/2\rfloor
 =\ell d-\left\lceil\frac{(k-2)\ell}{2}\right\rceil
 \le\ell d-\lceil\ell/2\rceil
\]
vertices of \(V_0\). Every cycle vector in the span is
supported on the union of the basis cycles, proving both
assertions about coverage.

For sharpness, take a theta graph, formed by three internally
disjoint paths of length \(\ell\) with the same endpoints. For even \(\ell\), put
both endpoints opposite \(V_0\); for odd \(\ell\), they
belong to opposite classes. The theta graph has cycle rank
two and \(\lfloor3\ell/2\rfloor\) vertices in \(V_0\).
Attach \(d-2\) otherwise disjoint \(2\ell\)-cycles at a
single vertex opposite \(V_0\). Each adds one to the cycle rank
and \(\ell\) to the vertex count, giving
\(\lfloor3\ell/2\rfloor+(d-2)\ell=\ell d-\lceil\ell/2\rceil\).
The shared cut vertex prevents additional simple cycles.

For the final assertions, recall that the cycle-space dimension of a
connected graph is \(|E|-|V|+1\). Distinct simple cycles give distinct
nonzero vectors, so a space of dimension two contains at most three
of them. Three such vectors in dimension two, or four in dimension
three, span the whole space. Applying the bound with \(\ell=3\)
gives four and seven, respectively.
\end{proof}

\begin{lemma}\label{k-lem:support-covered}\label{k-lem:cyclic-small-support}\label{k-lem:exceptional-local-weights}
In the bad-pair setting above, every vertex of \(S\) belongs to an
edge of \(J\). For \(s=6,7\), every incidence component of \(J\)
is cyclic; for \(s=7\), every component contains a complete triangle.
The weight satisfies
\[
\begin{array}{c|cccc}
s&5&6&7&8\\\hline
\text{lower bound on }W&3&6&8&10.
\end{array}
\]
\end{lemma}
\begin{proof}
We prove coverage first, then the assertions about components, and
finally the four weight bounds. Throughout, a \emph{private} support
vertex means a vertex lying in exactly one of
\(N_{\High}(p),N_{\High}(q)\).

\begin{claim}\label{k-claim:support-covered}
Every vertex of \(S\) belongs to an edge of \(J\).
\end{claim}
\begin{claimproof}
Suppose \(v\in S\) belongs to no edge of \(J\). Saturation of
the missing triple \(pqv\) forces, say, \(pv\) good. Let
\(y,z\in N_{\High}(p)\) be its additional core vertices.
The two edges at \(p\) represent \(py,pz\). Since no edge of
\(J\) contains \(v\), the representatives of \(vy,vz\) must
be the two edges at \(q\), namely \(qvy,qvz\). The closure
representing \(yz\) avoids these four used edges and hence avoids
\(I\). Reassigning \(qvy,qvz\) to \(qy,qz\) produces a witness
of Berge \(K_4\) for \(pq\), a contradiction.

\end{claimproof}

\begin{claim}\label{k-claim:support-cyclic}
For \(s=6,7\), every incidence component of \(J\) is cyclic. For \(s=7\), every component contains a complete triangle.
\end{claim}
\begin{claimproof}
If the two low supports are disjoint, every forced triangle is
complete. If they meet only at \(w\), the triangle forced at
\(w\) is complete. A noncomplete triangle at another vertex
\(v\) can use an edge of \(I\) only to represent \(vw\);
its other two representatives give a path from \(v\) to \(w\)
in the incidence graph of \(J\). Thus every component contains
a complete triangle in these two situations. They cover all cases
with \(s=7\), and all cases with \(s=6\) except the one in which
both supports have size four and their intersection has size two.

For that remaining case write the common vertices as \(u,v\) and
the private vertices as \(a,b,c,d\). Each four-point support is
partitioned into the two pairs in its incident hyperedges. The common
vertices can be paired together in both partitions, in neither, or in
exactly one. Up to interchanging \(p,q\), these give exactly
\[
\begin{array}{ll}
\mathrm{(i)}&I=\{puv,pab,quv,qcd\},\\
\mathrm{(ii)}&I=\{pua,pvb,quc,qvd\},\\
\mathrm{(iii)}&I=\{pua,pvb,quv,qcd\}.
\end{array}
\]
\begin{enumerate}
\item In (i), the triangles at \(a,b,c,d\) are complete. A
noncomplete triangle at \(u\) or \(v\) uses the pair \(uv\)
from \(I\); its remaining representatives connect it in \(J\)
to a private vertex. Hence every component contains a complete triangle.
\item In (ii), the pair \(uv\) has codegree zero. Indeed, any
representative of \(uv\), together with the four edges of \(I\),
would form a witness of Berge \(K_4\) for \(pq\). Therefore the
triangles forced at \(u,v\) are complete. Any noncomplete triangle
at a private vertex connects it in \(J\) to \(u\) or \(v\).
Again every component contains a complete triangle.
\item In (iii), a private vertex either lies in a complete triangle
or is connected to \(u\) or \(v\). A noncomplete triangle forced
at \(u\) has core \(uav\): it uses \(quv\) and two distinct
\(J\)-edges representing \(ua,av\). Similarly, one forced at
\(v\) gives two distinct representatives of \(ub,bv\).
If a central triangle is complete, its component is cyclic, and a
noncomplete triangle at the other central vertex joins that component.
If neither is complete, the two incidence paths from \(u\) to
\(v\), one through \(a\) and the other through \(b\), are
distinct and their union contains a cycle. This remains true if they
share edge nodes, since their internal support vertices differ.
Every other vertex either joins a central component or lies in a
complete triangle. Thus every component is cyclic.
\end{enumerate}

\end{claimproof}

\begin{claim}\label{k-claim:weight-five-six}
If \(s=5\), then \(W\ge3\); if \(s=6\), then \(W\ge6\).
\end{claim}
\begin{claimproof}
Edges of total weight at most two cover at most four vertices, so
coverage gives \(W\ge3\) when \(s=5\).
For general \(s\), let \(c\) be the number of incidence components.
The incidence graph has \(2h_2+3h_3\) edges and
\(s+h_2+h_3\) vertices. Its total cycle rank is consequently
\[
 (2h_2+3h_3)-(s+h_2+h_3)+c=W-s+c.
\]
When every component is cyclic, its cycle rank is at least one.
Claim~\ref{k-claim:support-cyclic}
therefore gives \(W\ge s\) for \(s=6,7\), proving the claimed
bound for \(s=6\).

\end{claimproof}

\begin{claim}\label{k-claim:weight-seven}
If \(s=7\), then \(W\ge8\).
\end{claim}
\begin{claimproof}
Since every component is cyclic, \(W\ge s=7\).
Suppose \(W=7\). The total cycle rank then equals the number of
incidence components, so every component is unicyclic. A unicyclic component has only one
simple cycle, so two distinct complete triangle cores must lie in
different components and in particular must be disjoint.
If the low supports are disjoint, all seven support vertices lie
on complete triangles. Disjoint three-point cores cannot cover seven
vertices. Otherwise, after relabelling, \(I=\{pwa,pbc,qwd,qef\}\). The missing triples \(pqw,pef,qbc\) force complete triangle cores
\(T_w,T_p,T_q\), respectively. For the latter two triples, the
pairs \(ef,bc\) are bad because they support type~2 edges.
The possible vertices of the forced cores are as follows:
\[
\begin{array}{c|l}
T_w&\{w,a,\xi\},\ \xi\in\{b,c\},\quad\text{or}
       \quad\{w,d,\eta\},\ \eta\in\{e,f\},\\
T_p&\text{one vertex from each of }\{e,f\},\{w,a\},\{b,c\},\\
T_q&\text{one vertex from each of }\{b,c\},\{w,d\},\{e,f\}.
\end{array}
\]
These triangles are complete: both edges at the selected low endpoint
are already used, and neither support pair of the other low vertex
is contained in the corresponding displayed core.
The core \(T_w\) differs from both \(T_p,T_q\).
If \(T_p=T_q\), their common core contains \(w\), contradicting
its required disjointness from \(T_w\). If \(T_p\ne T_q\), both
must avoid \(w\), so \(a\in T_p\) and \(d\in T_q\).
Then either form of \(T_w\) intersects one of them, again a
contradiction. Hence \(W\ge8\).

\end{claimproof}

\begin{claim}\label{k-claim:weight-eight}
If \(s=8\), then \(W\ge10\).
\end{claim}
\begin{claimproof}
Write the four edges of \(I\) as
\(\{p\}\cup P_1,\{p\}\cup P_2,\{q\}\cup Q_1,\{q\}\cup Q_2\),
where the four pairs \(P_1,P_2,Q_1,Q_2\) are pairwise disjoint.
Every support vertex lies in a complete triangle. In addition, the
four missing triples \(\{p\}\cup Q_i,\{q\}\cup P_i\), \(i=1,2\),
force four distinct complete triangle cores. The cores forced by \(\{p\}\cup Q_i\) contain one vertex of
\(Q_i\) and one from each of \(P_1,P_2\); those forced by
\(\{q\}\cup P_i\) contain one vertex of \(P_i\) and one from
each of \(Q_1,Q_2\).
These descriptions establish both completeness and distinctness.

Every component contains at least three support vertices, so there
are at most two components. All are cyclic, and hence \(W\ge8\).
If \(W=8\), every component is unicyclic. Since every support
vertex lies on a six-cycle, each such component has exactly three
support vertices, which cannot account for eight vertices.
If \(W=9\) and the incidence graph is connected, its cycle rank is two,
but it has at least four distinct six-cycles, impossible by
Theorem~\ref{k-thm:short-cycle-support}. With two components, their
cycle ranks sum to three and hence are one and two. The component of cycle rank one
has exactly three support vertices and one complete triangle core. The other
has five support vertices, all lying on cycles, and at least three
distinct six-cycles. The same theorem allows at most four support
vertices on its cycles. This final contradiction proves \(W\ge10\).
\end{claimproof}

The five claims prove all the assertions of the lemma.
\end{proof}

\section{Proof of the Berge-\texorpdfstring{\(K_4\)}{K4} saturation formula}\label{k-sec:constructions}
We now combine the local reductions with the support estimates in
one proof of the exact formula.

\begin{proof}[Proof of Theorem~\ref{k-thm:main}]
The upper bounds follow from Chen, He, and Xie~\cite{CHX}.
Their Constructions~2.1 and~2.2, together with their Lemma~2.1,
give \(\sat_3(n,\KF)\le n\) for every \(n\ge7\).
Their tight five-cycle, the 3-graph on \(\{0,1,2,3,4\}\) with
edges \(012,123,234,034,014\), and its union with an isolated
vertex give the upper bound five for \(n=5,6\); these exact values are included
in~\cite[Theorem~1.1]{CHX}.

We begin the lower bound with the small orders.
For \(n\le4\), even the complete 3-graph has fewer than six edges
and is \(\KF\)-free; hence saturation requires all \(\binom n3\) triples.
For \(n=5,6\), fewer than five edges cannot be saturated, since
adding a missing triple still leaves fewer than six edges.
The case \(n=8\) is Proposition~\ref{k-prop:typeone-orders-eight-nine}.

If the bound fails at \(n=7\), choose a minimum-edge saturated 3-graph.
Then \(m\le6\), \(\lambda\ge3\), and \(x\le4\).
Proposition~\ref{k-prop:typeone} makes the type~1 graph a matching,
and Lemma~\ref{k-lem:twohigh} gives \(q_{\High}\ge2\).
The matching count yields \(m\ge3\lambda/2+2\), so
\(7=\lambda+x\le3x-6\le6\), a contradiction.

It remains to consider \(n\ge9\). Suppose the lower bound fails,
and choose an \(n\)-vertex saturated 3-graph with
\(m=\sat_3(n,\KF)\le n-1\). This is a counterexample in the
sense of Section~\ref{k-sec:lower}; in particular, \(\lambda\ge3\).
Proposition~\ref{k-cor:no-type-one} excludes type~1 edges.
We first force a bad pair of low vertices, then eliminate every
possible size of its support.

\begin{claim}\label{k-prop:allgood}
There is a bad pair in \(\Low\).
\end{claim}
\begin{claimproof}
Suppose every pair in \(\Low\) is good. Fix \(a_0\in\Low\)
with support pairs \(L,R\). For each \(b\in \Low\setminus\{a_0\}\),
choose a witness of Berge \(K_4\) for \(a_0b\), with additional pair \(P_b\)
meeting both \(L,R\), and closure \(f_b\supseteq P_b\)
avoiding \(a_0,b\). Both edges at \(b\) meet \(P_b\).
Put \(D=\bigcup_bP_b\), \(d=|D|\in\{2,3,4\}\), and
\(t_0=|L\cap D|+|R\cap D|\ge d\). Write \(\sigma_D=\sum_{v\in D}d(v)\). Incidence counting gives
\(\sigma_D=2(\lambda-1)+t_0+C\).
Here \(C\) counts the extra \(D\)-incidences in edges at
\(\Low\setminus\{a_0\}\), beyond one per such edge, together with
all \(D\)-incidences in all-high edges. The edges at \(a_0\) are
excluded from \(C\), since all their incidences are already counted
in \(t_0\). Every summand is nonnegative: each edge at
\(b\in\Low\setminus\{a_0\}\) meets \(P_b\subseteq D\).
We have \(C\ge2\): an all-high closure contributes at least two;
otherwise \(f_b=cP_b\) contributes one and the distinct
closure for \(a_0c\) contributes another.

For every \(d\), this incidence identity and degree summation give \(3m\ge2\lambda+3(x-d)+\sigma_D =3n+\lambda+t_0+C-3d-2\).

\smallskip\noindent\emph{Case \(d=2\).}
The inequality gives \(3m\ge3n+\lambda-4>3n-3\), since \(\lambda\ge3\),
contradicting \(m\le n-1\).

\smallskip\noindent\emph{Case \(d=3\).}
The same inequality gives
\(3m\ge3n+\lambda-6\), which is impossible unless
\(\lambda=3,t_0=3\). In that case all high vertices are cubic:
the inequalities \(3m\le3n-3=2\lambda+3x\le3m\) force equality
in every high-vertex degree bound.
Write the edges at \(a_0\) as \(a_0pq,a_0rs\), with
\(D=\{p,r,s\}\), and choose \(P_b=pr,P_c=ps\).
The same incidence identity forces \(C=2\). An all-high closure containing all of \(D\) contributes three
incidences. Otherwise it contains only its selected pair in \(D\)
and cannot close the other selected pair, so the other closure
contributes at least one further incidence. Both alternatives
contradict \(C=2\). Hence the closures are exactly \(cpr,bps\). The other edges at \(b\) and \(c\) are
\(bru,csv\), with \(u,v\notin D\). The pairs \(pq,rs\)
occur only in the respective edges at \(a_0\). A witness of
Berge \(K_4\) for \(a_0p\) would use \(a_0pq\) for \(a_0q\), leaving
\(pq\) unrepresented; similarly \(a_0r\) is bad.
The support \(pr\) of the type~2 edge \(cpr\) is bad by
Lemma~\ref{k-lem:goodfacts}(1), so the missing triple \(a_0pr\) contains no good pair, contradicting saturation.

\smallskip\noindent\emph{Case \(d=4\).}
The supports \(L,R\) are disjoint. The selected pairs each meet both \(L,R\) and together cover all
four vertices. Since both parts have size two, two selected pairs,
say \(P_b,P_c\), are disjoint.
Their closures are distinct, and \(C\ge4\). Indeed, two
all-high closures suffice. If \(f_c=zP_c\) is a type~2 edge, then
\(z\notin\{a_0,b,c\}\). In particular, \(z\ne b\), because
every edge at \(b\) meets \(P_b\), whereas \(P_c\cap P_b=\varnothing\).
Goodness of \(bz\) forces an edge
at \(b\) to meet both \(P_b,P_c\), contributing an extra
incidence outside the closures. If both closures are type~2 edges,
the symmetric argument gives a further edge at \(c\).
Thus \(\sigma_D\ge2\lambda+6\), and \(3m\ge3n+\lambda-6\), excluding
\(\lambda\ge4\).

If \(\lambda=3\), all high vertices are cubic. Write
\(P_b=pr,P_c=qs\). Both closures are all-high, and
\(C=4\) forces the relevant edges to be
\(
 a_0pq,a_0rs,bpu,brv,cqw,csz,prt,qsh,
\)
where \(u,v,w,z,t,h\notin D\).
The two additional core vertices in a witness of Berge \(K_4\)
for \(bc\) must belong to \(\{u,v\}\cap\{w,z\}\), since
\(\{p,r\}\cap\{q,s\}=\varnothing\) and the displayed outside
vertices avoid \(D\). Thus goodness of \(bc\) forces
\(\{u,v\}=\{w,z\}\) and \(u\ne v\). The common high neighbours of \(p,r\),
apart from the endpoints, are then contained in \(\{t\}\),
so \(pr\) is bad. As above, the unique representatives of
\(pq,rs\) make \(a_0p,a_0r\) bad. Missing \(a_0pr\)
again contradicts saturation.
\end{claimproof}

Choose a bad pair \(p,q\in\Low\), and form \(S,I,J,W\) as in
Subsection~\ref{k-sec:bad-low}. The remaining argument is organized
by \(s=|S|\). We will repeatedly use
\[
 s\le8,\quad x\ge\max\{6,s\},\quad n\le2x-3,
 \quad x\le4s-9-\lceil3W/2\rceil.
\]
The last inequality is the incidence bound established in Subsection~\ref{k-sec:bad-low}.
For \(s\ge5\), combining it with the weight bounds from
Lemma~\ref{k-lem:exceptional-local-weights} leaves only the following
parameter values:
\[
\begin{array}{c|c|c|c}
s&\text{lower bound on }W&\text{upper bound on }x&
\text{remaining }(x,W)\\\hline
5&3&6&(6,3)\\
6&6&6&(6,6)\\
7&8&7&(7,8)\\
8&10&8&(8,10).
\end{array}
\]
For example, if \(s=6\), then \(x\le15-\lceil3W/2\rceil\le6\).
Since \(x\ge6\), increasing \(W\) above six is impossible.
The other rows follow in the same way.
Claim~\ref{k-lem:small-support-sharp} also excludes \(s\le4\);
the subsequent support claims address the three larger sizes.

We record one component observation that will be used in both the
\(s=7\) and \(s=8\) cases.
\begin{claim}\label{k-claim:component-parity}
If \(s\in\{7,8\}\) and \(S=\High\), then every incidence
component of \(J\) contains an all-high edge and has even weight
at least four.
\end{claim}
\begin{claimproof}
For \(s=7\), every component contains a complete triangle by
Lemma~\ref{k-lem:support-covered}; for \(s=8\), every support
vertex lies in one. Suppose a component has no all-high edge.
Its complete triangle is then represented by three type~2 edges,
so all three core pairs are bad by Lemma~\ref{k-lem:goodfacts}(1).
The core triple cannot be an edge: such an all-high edge would also
belong to this component. The missing core triple thus contradicts
saturation.

It remains to prove the parity assertion. Because \(S=\High\),
every pair projection in \(J\) comes from an edge at a vertex
\(r\in\Low\setminus\{p,q\}\). The missing triple \(pqr\)
forces \(pr\) or \(qr\) good, say \(pr\). Its additional core
vertices lie in \(N_{\High}(p)\) and are supplied by the two
distinct edges at \(r\). Its closure avoids the used edges at
\(p,r\). It also avoids the edges at \(q\), because the two
low supports intersect in at most one vertex, whereas a closure
contains both additional core vertices. The closure therefore gives an
edge of \(J\) joining the two projections at \(r\) into the same
component. Thus pair projections occur in pairs within each component.
Triple projections have weight two, so every component has even weight.
Since it contains a complete triangle, its weight is at least three
and hence at least four.
\end{claimproof}

\begin{claim}\label{k-lem:small-support-sharp}\label{k-lem:support-five-impossible}
The support has \(s\ge6\) vertices.
\end{claim}
\begin{claimproof}
For any triangle forced at a support vertex, the two edges at its
selected low endpoint are already used. Only two edges of \(I\)
remain, whereas the triangle needs three representatives. At least
one representative therefore avoids \(I\), and \(W\ge1\).
If \(s\le4\), the bound \(x\le4s-9-\lceil3W/2\rceil\) gives
\(x\le16-9-2=5\), contrary to \(x\ge6\).

Suppose \(s=5\). We divide the argument into the parameter reduction,
the shape of \(J\), and the final obstruction.

\smallskip\noindent\emph{Forced parameters.}
The parameter table gives \(x=6,W=3\), and
\(n\le2x-3=9\). Since \(n\ge9\), we have \(n=9\) and \(\lambda=3\).
Now \(m=2\lambda+q_{\High}\), \(q_{\High}\ge2\), and \(m\le n-1\)
give \(q_{\High}=2,m=8\). Moreover,
\(3m=24=2\lambda+3x\), so every high vertex is cubic.
Write \(\Low=\{p,q,r\}\) and \(\High\setminus S=\{z\}\).
The four edges outside \(I\) are the two edges at \(r\) and
the two all-high edges. Exactly three contain \(z\).

If exactly one all-high edge contains \(z\), both edges at \(r\)
contain it. The latter edges meet \(S\) only once and do not occur
in \(J\). Thus \(J\) consists of a pair and a triple, from the
two all-high edges. Coverage of all five support vertices makes
these projections disjoint. The original all-high edges are then
also disjoint, contradicting Lemma~\ref{k-lem:two-high-not-disjoint}.

\smallskip\noindent\emph{The three projected pairs.}
Both all-high edges must therefore contain \(z\), and exactly one
edge at \(r\) contains it. Write the edges at \(r\) as
\(rzc\) and \(rL\), where \(c\in S\) and \(L\) is a pair in
\(S\). The two all-high edges and \(rL\) give three pair edges
in \(J\). For \(v\in S\), put
\(\gamma(v)=|\{e\in I:v\in e\}|\). Cubicity gives \(d_J(c)=2-\gamma(c)\) and
\(d_J(v)=3-\gamma(v)\) for \(v\ne c\). Each \(\gamma(v)\ge1\), and coverage gives \(d_J(v)\ge1\).
Consequently \(d_J(c)=1\), and every other \(J\)-degree is one
or two. The total degree is six on five vertices, so exactly one
vertex has degree two. Hence \(J\) is a two-edge path together
with an isolated pair. Parallel pairs are impossible because they
would give two vertices degree at least two.

The missing triple \(pqr\) forces, after interchanging \(p,q\)
if needed, \(pr\) good. Its additional core vertices are \(c\) and
a vertex of \(L\): they must lie in \(S\), and the two edges
at \(r\) must supply them separately. Since \(d_J(c)=1\),
we have \(\gamma(c)=1\). This one incidence belongs to an edge
at \(p\), so no edge at \(q\) contains \(c\). The closure
therefore avoids \(I\); it also avoids both used edges at \(r\),
so it is all-high. Its projection joins \(c\) to an endpoint of
\(L\). Thus \(L\) and the closure projection form the two-edge
path, with \(c\) as an endpoint. Relabel so that the four edges
outside \(I\) are \(rzu,\; rvw,\; zuv,\; zab\). Here \(S=\{u,v,w,a,b\}\). The vertices \(u,v\) each occur once
in \(I\), and \(w,a,b\) each occur twice.

\smallskip\noindent\emph{The incident edges and a bad triple.}
The two edges at \(p\) separately represent \(pu,pv\), so they
are \(pu\alpha,pv\beta\), with
\(\alpha,\beta\in\{w,a,b\}\). If \(\alpha=\beta\), the
other two vertices of \(\{w,a,b\}\) would both have to occur
in both edges at \(q\), repeating a hyperedge. Hence
\(\alpha\ne\beta\), and the remaining incidences force \(I=\{pu\alpha,pv\beta,qt\alpha,qt\beta\}, \; \{t,\alpha,\beta\}=\{w,a,b\}\). The three pairs \(qt,q\alpha,q\beta\) are bad. Indeed, a
witness of Berge \(K_4\) for any one of them must use both
edges at \(q\) and would then need two additional representatives
containing \(t\). Cubicity of \(t\) leaves only one unused edge
at \(t\). Saturation of the missing triple \(q\alpha\beta\)
therefore forces \(\alpha\beta\) good.

If \(t=w\), then \(\{\alpha,\beta\}=\{a,b\}\), and the five
distinct edges \(pu\alpha,\; pv\beta,\; qw\alpha,\; qw\beta,\; zab\) represent \(p\alpha,p\beta,q\alpha,q\beta,\alpha\beta\).
They form a witness of Berge \(K_4\) for \(pq\), contrary to its
choice. Otherwise we may take \(t=a\), since \(t=b\) is symmetric.
The pair forced to be good is then \(wb\). Its common high neighbours
are determined by the eight edges specified above:
\[
\begin{array}{c|c|l}
(\alpha,\beta)&N_{\High}(w)\cap N_{\High}(b)&\text{obstruction}\\\hline
(w,b)&\{v,a\}&\text{no hyperedge contains }va,\\
(b,w)&\{a\}&\text{only one possible additional core vertex}.
\end{array}
\]
In the first row a closure on \(va\) is impossible; in the second
there are not two additional core vertices. Thus \(wb\) is bad
in both cases, completing the contradiction.
\end{claimproof}

\begin{claim}\label{k-lem:support-six-impossible}
The case \(s=6\) is impossible.
\end{claim}
\begin{claimproof}
The parameter table gives \(x=6,W=6\), so \(\High=S\).
All edges outside \(I\) now project to either pairs at the other
low vertices or all-high triples. Hence \(6=W=2(\lambda-2)+2q_{\High}\). Since \(\lambda\ge3\) and \(q_{\High}\ge2\), equality forces
\(\lambda=3,q_{\High}=2\), and \(m=2\lambda+q_{\High}=8\).
All high vertices are cubic, and \(J\) has two pair edges and
two triple edges. For \(v\in S\), write
\(\gamma(v)=|\{e\in I:v\in e\}|\); then \(d_J(v)=3-\gamma(v)\). Every incidence component is cyclic by Lemma~\ref{k-lem:support-covered}.
Since \(W=s=6\), the total cycle rank equals the number of
incidence components. Every component is cyclic, so every component
is unicyclic.
In particular, a component has at most one complete triangle core.

Let \(A_p=N_{\High}(p)\), \(A_q=N_{\High}(q)\).
Since each support has size three or four, the possibilities for
\((|A_p|,|A_q|,|A_p\cap A_q|)\), up to interchanging \(p,q\), are \((3,3,0),\; (3,4,1),\; (4,4,2)\). \begin{enumerate}
\item If the supports are disjoint, every support vertex belongs
to a complete triangle. In a three-point support, the two support
pairs have a common vertex. That vertex has \(\gamma=2\) and
\(d_J=1\), so it cannot lie on a six-cycle, a contradiction.
\item If the supports meet in one vertex, the triangle forced at
that vertex is complete. But it belongs to at least two edges of
\(I\), so again \(d_J\le1\), a contradiction.
\item Suppose both supports have four vertices and share \(u,v\).
Use partitions (i)--(iii) in the proof of Claim~\ref{k-claim:support-cyclic}.
In partition (ii), \(u,v\) have \(d_J=1\) but their forced
triangles are complete. In partition (i), all four private vertices
lie in complete triangles, so at least two distinct complete triangle cores
are needed. They cannot lie in the same unicyclic component.
Two different components containing such triangles require at least
six edge nodes, whereas \(J\) has only four.

It remains to consider partition (iii), namely
\(I=\{pua,pvb,quv,qcd\}\). Here
\(d_J(u)=d_J(v)=1\) and
\(d_J(a)=d_J(b)=d_J(c)=d_J(d)=2\).
The forced triangles at \(u,v\) are noncomplete. Their forms,
established in the proof of Claim~\ref{k-claim:support-cyclic}, require
representatives of \(ua,av\) and of \(ub,bv\), respectively.
The sole projected edge at \(u\) must contain both \(ua\) and
\(ub\), so it is \(uab\). Similarly, the sole edge at \(v\)
is \(vab\). These two triple projections exhaust the degrees of
\(a,b\). The two remaining pair projections must both be \(cd\).
Both come from the only other low vertex \(r\), and would require
the repeated hyperedge \(rcd\), impossible in a simple 3-graph.
\end{enumerate}
\end{claimproof}

\Needspace{6\baselineskip}
\begin{claim}\label{k-lem:support-seven-new}
The case \(s=7\) is impossible.
\end{claim}
\begin{claimproof}
The parameter table gives \(x=7,W=8\), so \(\High=S\).
By Claim~\ref{k-claim:component-parity}, every incidence component
contains an all-high edge and has even weight at least four. Thus
there are either one or two components.

\smallskip\noindent\emph{At least three complete triangle cores are required.}
If the low supports are disjoint, all seven support vertices lie in
complete triangles. At least three distinct cores are needed to
cover them. If the supports intersect, their sizes are both four
and their intersection has one vertex. Write \(I=\{pwa,pbc,qwd,qef\}\). A noncomplete triangle can borrow an \(I\)-edge only along
\(wa\) or \(wd\), and it can be forced only at \(a\) or
\(d\). Thus each of \(w,b,c,e,f\) lies in a complete triangle.
Use the cores \(T_w,T_p,T_q\) forced by \(pqw,pef,qbc\), as
listed in the proof of Claim~\ref{k-claim:weight-seven}.
The core \(T_w\) differs from both of the others. If
\(T_p\ne T_q\), there are already three distinct cores.
If \(T_p=T_q\), the common core consists of \(w\), one vertex
of \(\{b,c\}\), and one of \(\{e,f\}\). The union
\(T_w\cup T_p\) contains at most three of \(b,c,e,f\).
A missed vertex belongs to another complete triangle, giving a
third core. In both support configurations at least five support
vertices lie on complete triangles.

\smallskip\noindent\emph{Neither component count is possible.}
If there are two components, both have weight four. Each contains
an all-high edge and a complete triangle. It cannot consist of two
triple edges, since a Berge triangle needs three distinct edge nodes.
It therefore consists of one triple and two pairs. For these three
edges to represent a triangle, the pairs must be distinct and share
one endpoint; they determine its unique core. The two components
can consequently contain at most two complete triangle cores, contrary to the preceding conclusion.

If the incidence graph is connected, its cycle rank is
\(W-s+1=2\). The preceding argument gives three distinct six-cycles and at least
five support vertices lying on cycles. These three distinct cycle vectors necessarily span the
two-dimensional cycle space. Theorem~\ref{k-thm:short-cycle-support}
therefore shows that at most four support vertices lie on cycles. This contradiction excludes \(s=7\).
\end{claimproof}

\begin{claim}\label{k-lem:eight}
The case \(s=8\) is impossible.
\end{claim}
\begin{claimproof}
The parameter table gives \(x=8,W=10\), and hence \(\High=S\).
Every support vertex lies on a six-cycle, and the four missing
triples used in Claim~\ref{k-claim:weight-eight} supply
four distinct complete triangle cores.

If the incidence graph were connected, its cycle rank would be
\(W-s+1=3\). Four distinct nonzero cycle vectors cannot lie in
a space of dimension two, which has only three nonzero vectors.
The four six-cycles therefore span the full cycle space.
Theorem~\ref{k-thm:short-cycle-support} shows that at most seven
support vertices lie on cycles, contrary to the coverage of all eight.

Thus there are at least two components. Each contains at least three
support vertices, so there are exactly two. By
Claim~\ref{k-claim:component-parity}, both have even weight at least
four, and their weights sum to ten. Their weights are therefore
four and six.

Consider the weight-four component. It contains an all-high edge
and a complete triangle, so, as in the preceding claim, it has one
triple and two adjacent pair edges. Write these pairs as \(ab,bc\).
The triple must contain \(ac\), in order to complete the triangle.
If its third vertex were outside \(\{a,b,c\}\), that vertex would
have incidence degree one and lie on no cycle. Since every support
vertex lies on a six-cycle, the triple is \(abc\).
This component has exactly three support vertices and one complete triangle core.

The other component therefore has five support vertices and weight
six, giving cycle rank \(6-5+1=2\). At least three of the four
distinct complete triangle cores belong to it, and all five of its support
vertices lie on cycles. The four-vertex bound in
Theorem~\ref{k-thm:short-cycle-support} gives the final contradiction.
\end{claimproof}

The support cases exhaust all possibilities for a counterexample.
Together with the small-order argument and the cited constructions,
this proves the asserted saturation formula.
\end{proof}

\begin{remark}
The five-star result determines both endpoints and every intermediate
edge count, including the exceptional gaps immediately below the
Tur\'an number. Its proof separates a structural restriction on dense
saturated hypergraphs from a finite construction problem. For Berge-\(K_4\),
the incidence inequalities and support analysis establish the lower bound
\(m\ge n\) throughout its exact range of validity. The known upper-bound
constructions then yield the saturation formula at every order.
These two arguments use distinct local completion certificates:
witnesses of Berge five-star in the star problem and witnesses of
Berge \(K_4\) in the clique problem. Together with the structural
counts, they give the stated exact saturation results.
\end{remark}

\end{document}